\documentclass[11pt,reqno,a4paper]{amsart}
\usepackage{color}
\usepackage{amssymb,amsmath}
\usepackage[latin1]{inputenc}
\usepackage[active]{srcltx}
\usepackage{graphicx}
\usepackage{exscale,relsize}
\usepackage{textgreek}
\usepackage{epsfig,graphics}
\usepackage{psfrag}
\usepackage{caption}
\usepackage{subcaption}
\usepackage[toc,page]{appendix}
\usepackage{bigints}
\def\N{\mathbb{N}}
\def\R{\mathbb{R}}
\def\m1{{I\!\!M}}

\newcommand{\grad}{\nabla}

\renewcommand{\to}{\rightarrow}
\newcommand{\pa}{\partial}

\newcommand{\ino}{\int_{\Omega}}

\newcommand{\ainf}{\mbox{as\;}\;n\to+\infty}

\newcommand{\fo}{\forall}

\newcommand{\rife}[1]{(\ref{#1})}
\newcommand{\ov}[1]{\overline{#1}}
\newcommand{\un}[1]{\underline{#1}}

\newcommand{\scp}{\scriptstyle}
\newcommand{\sscp}{\scriptscriptstyle}
\newcommand{\dsp}{\displaystyle}

\renewcommand{\dfrac}{\displaystyle\frac}
\newcommand{\finedim}{\hspace{\fill}$\square$}
\newcommand{\intbar}{\mathop{\int\makebox(-15.5,0){\rule[6pt]{.7em}{0.3pt}}\kern-6pt}\nolimits}

\newcommand{\ii}{\infty}

\newcommand{\eps}{\varepsilon}
\newcommand{\dt}{\delta}

\newcommand{\al}{\alpha}
\newcommand{\be}{\beta}

\newcommand{\wsg}{\widehat{\sigma}}
\newcommand{\sg}{\sigma}
\newcommand{\ga}{\gamma}
\newcommand{\om}{\Omega}
\newcommand{\lm}{\lambda}

\newcommand{\wlm}{\widehat{\lm}_1}
\newcommand{\wpl}{\widehat{\psi}}
\newcommand{\wall}{\widehat{\al}}

\newcommand{\rl}{\mbox{\Large \textrho}_{\!\sscp \lm}}
\newcommand{\rla}{\mbox{\Large \textrho}_{\!\sscp \lm,\sscp \al}}

\newcommand{\rlq}{(\mbox{\Large \textrho}_{\! \sscp \lm})^{\frac1q}}

\newcommand{\rlqa}{(\mbox{\Large \textrho}_{\! \sscp \lm,\al})^{\frac1q}}

\renewcommand{\rho}{\mbox{\Large \textrho}}
\newcommand{\rh}{\mbox{\Large \textrho}}

\newcommand{\pl}{\psi_{\sscp \lm}}

\newcommand{\um}{u_{\sscp \mu}}
\newcommand{\val}{v_{\sscp I}}

\newcommand{\vxi}{\xi}

\newcommand{\vl}{v_{\sscp \lm}}
\newcommand{\bl}{\be_{\sscp \lm}}

\newcommand{\ssl}{\sscp \lm}
\newcommand{\ml}{m_{\sscp \lm}}
\newcommand{\all}{\al_{\ssl}}
\newcommand{\gal}{\ga_{\sscp I}}
\newcommand{\el}{E_{\ssl}}

\newcommand{\fbi}{{\bf (F)$_{I}$}}
\newcommand{\prl}{{\textbf{(}\mathbf P\textbf{)}_{\mathbf \lm}}}

\newtheorem{theorem}{Theorem}[section]
\newtheorem{proposition}[theorem]{Proposition}
\newtheorem{lemma}[theorem]{Lemma}
\newtheorem{corollary}[theorem]{Corollary}
\newtheorem{remark}[theorem]{Remark}
\newtheorem{definition}[theorem]{Definition}
\newcommand{\brm}{\begin{remark}\rm}
\newcommand{\erm}{\end{remark}}
\newcommand{\bdf}{\begin{definition}\rm}
\newcommand{\edf}{\end{definition}}
\newcommand{\bte}{\begin{theorem}}
\newcommand{\ete}{\end{theorem}}
\newcommand{\bpr}{\begin{proposition}}
\newcommand{\epr}{\end{proposition}}
\newcommand{\ble}{\begin{lemma}}
\newcommand{\ele}{\end{lemma}}
\newcommand{\bco}{\begin{corollary}}
\newcommand{\eco}{\end{corollary}}
\newcommand{\beq}{\begin{equation}}
\newcommand{\eeq}{\end{equation}}
\newcommand{\bdm}{\begin{displaymath}}
\newcommand{\edm}{\end{displaymath}}

\newcommand{\graf}[1]{\left\{\begin{array}{ll}#1\end{array}\right.}

\def\sideremark#1{\ifvmode\leavevmode\fi\vadjust{\vbox to0pt{\vss
 \hbox to 0pt{\hskip\hsize\hskip1em \vbox{\hsize2.1cm\tiny\raggedright\pretolerance10000 \noindent #1\hfill}\hss}\vbox to15pt{\vfil}\vss}}}

\begin{document}
\numberwithin{equation}{section}
\parindent=0pt
\hfuzz=2pt
\frenchspacing

\title[Generic properties for free boundary problems]{Generic properties of free boundary problems \\in plasma physics}

\author[D. Bartolucci]{Daniele Bartolucci}
\address{Daniele Bartolucci, Department of Mathematics, University of Rome ``Tor Vergata", Via della ricerca scientifica 1, 00133 Roma, Italy}
\email{bartoluc@mat.uniroma2.it}

\author[Y. Hu]{Yeyao Hu}
\address{Yeyao Hu, School of Mathematics and Statistics, HNP-LAMA, Central South University,
Changsha, Hunan 410083, P.R. China}
\email{huyeyao@csu.edu.cn}

\author[A. Jevnikar]{ Aleks Jevnikar}
\address{Aleks Jevnikar, Department of Mathematics, Computer Science and Physics, University of Udine, Via delle Scienze 206, 33100 Udine, Italy}
\email{aleks.jevnikar@uniud.it}

\author[W. Yang]{Wen Yang}
\address{Wen Yang, Wuhan Institute of Physics and Mathematics, Chinese Academy of Sciences, P.O. Box 71010,
Wuhan 430071, P.R. China}
\address{Innovation Academy for Precision Measurement Science and Technology, Chinese Academy of
Sciences, Wuhan 430071, P.R. China}
\email{wyang@wipm.ac.cn}


\thanks{D.B. has been partially supported by Beyond Borders project 2019 (sponsored by Univ. of Rome "Tor Vergata") "\emph{Variational Approaches to PDE's}", MIUR Excellence Department Project awarded to the Department of Mathematics, Univ. of Rome Tor Vergata, CUP E83C18000100006.}

\begin{abstract}
We are concerned with the global bifurcation analysis of positive solutions to free boundary problems arising in plasma physics. We show that in general, in the sense of domain variations, the following alternative holds: either the shape of the branch of solutions resembles the monotone one of the model case of the two-dimensional disk, or it is a continuous simple curve without bifurcation points which ends up at a point where the boundary density vanishes. On the other hand, we deduce a general criterion ensuring the existence of a free boundary in the interior of the domain. Application to a classic nonlinear eigenvalue problem is also discussed.
\end{abstract}
\maketitle

\tableofcontents


\setcounter{section}{0}
\section{\bf Introduction}
\setcounter{equation}{0}

\bigskip
\bigskip

For {$\om\subset \R^N$}, $N\geq2$, an open and bounded domain of class $C^{4}$ we consider the free boundary problem
$$
\graf{-\Delta v = (v)_{+}^p\quad \mbox{in}\;\;\om\\ \\
-\bigints\limits_{\pa\om} \dfrac{\pa v}{\pa\nu}=I \\ \\
v=\ga \quad \mbox{on}\;\;\pa\om
}\qquad \qquad \mbox{\bf (F)}_{I}
$$
for the unknowns $\ga \in \R$ and $v\in C^{2,r}(\ov{\om}\,)$, $r\in (0,1)$.
Here $(v)_+$ is the positive part of $v$, $\nu$ is the exterior unit normal, $I>0$ and $p\in (1,p_{N})$ are fixed, with
$$
p_{ N}=\graf{+\ii,\;N=2 \vspace{0.2cm}\\ \frac{N}{N-2}\,,\; N\geq 3.}
$$
If not empty, the set $\om \cap \partial\{x\in\Omega \,|\, v(x)>0\}$ is called the free boundary of the problem. For $p=1$ one needs a slightly different formulation which we skip here since such model case is by now fully understood. Up to a suitable scaling, one can assume without loss of generality
that $|\om|=1$. In particular, we could {pick} any positive constant to multiply $(v)_+^p$.

\

Problems like \fbi\, arise in modeling the plasma physics of Tokamak reactors \cite{Fre,Kad,Mer,Wes}. The rigorous
mathematical analysis was first undertaken in \cite{BeBr,Te,Te2}, see also \cite{BJ1} for a complete list of references.

\subsection{The Grad-Shafranov plasma equilibrium equation in a Tokamak}
In this subsection we first sketch the derivation of the Grad-Shafranov equation following  \cite{blum} and \cite{Te}, see \cite{Fre,Wes} for more details and references. Then we will discuss its
connection with problem \fbi\, and the physical motivations behind our results.\\

In fusion experiments a plasma of ions of deuterium and tritium and their electrons is confined to
a toroidal region by a magnetic field. The vessel containing the plasma is itself a toroid and the magnetic field must be strong enough to guarantee that the high temperature (about $10$ KeV which is about $10$ x $ 10^7$ centigrade degree) plasma does not touch the vessel wall.\\
Let ${\bf j}$ be the plasma current density, {\bf B} the magnetic field and $p$ the kinetic pressure. Then the condition for equilibrium in the region occupied by the plasma requires that the magnetic force balances the force due to the plasma pressure,
\beq\label{equi1}
{\bf j}\times {\bf B}=\grad p,
\eeq
to be solved together with Ampere's Theorem,
\beq\label{equi2}
\mbox{curl } {\bf B}=\mu_0{\bf j},
\eeq
where $\mu_0$ is the {vacuum magnetic permeability}, and the conservation of magnetic induction,
\beq\label{equi3}
\mbox{div }{\bf B}=0.
\eeq
On the other side, in the region which is not occupied by the plasma we have to solve \rife{equi2} with ${\bf j}=0$
together with \rife{equi3}.\\

The Grad-Shafranov equation is a rewriting of \rife{equi1}, \rife{equi2}, \rife{equi3} in toroidal
coordinates under axisymmetric assumptions. The tokamak is a torus obtained by the rotation around the vertical axis
$\mbox{O}z$ of a two dimensional simply connected smooth and bounded domain $\om$ lying in the plane $\mbox{O}xz$ with $x$ positive. The plasma occupies
the toroidal region defined by the rotation of $\om_{+}\subseteq \om$ around the vertical axis
$\mbox{O}z$, and we set $\om_{-}=\{\om\setminus \om_{+}\} \subset \om$.\\
Consider the cylindrical coordinate system
$\{{\bf e}_r,{\bf e}_\varphi,{\bf e}_z\}$ where $0\leq \varphi< 2\pi$ and assume ${\bf B}={\bf B}(r,z)$ to be independent by the toroidal angle $\varphi$.
Then decompose ${\bf B}={\bf B}_{\perp}+{\bf B}_{\varphi}$, where
${\bf B}_{\perp}=B_r{\bf e}_r+B_z{\bf e}_z$ is the poloidal field and
${\bf B}_{\varphi}=B_{\varphi}{\bf e}_\varphi$ is the toroidal field. Next, introducing the (normalized) poloidal flux
$$
\psi(r,z)=\int\limits_0^rB_z(t,z) tdt,
$$
we have
\beq\label{bdy1}
B_z(r,z)=\frac{1}{r}\frac{\pa \psi}{\pa r}(r,z),
\eeq
and we deduce from \rife{equi3} that
\beq\label{bdy2}
B_r(r,z)=-\frac{1}{r}\frac{\pa \psi}{\pa z}(r,z).
\eeq
As a consequence, putting $\nabla =\left(\frac{\pa }{\pa r},0, \frac{\pa}{\pa z}\right)$, we have
${\bf B}_{\perp}=\frac{1}{r}\grad \psi \times {\bf e}_{\varphi}$ and in particular
${\bf B} \cdot \grad \psi=0$. Also, from \rife{equi1}, we see that
${\bf B} \cdot \grad p=0$. In other words the region inside the tokamak is foliated by a family of surfaces (called magnetic surfaces) and on each magnetic surface
$p$ and $\psi$ are constant, which is why one assumes that $p=p(\psi)$. Since from \rife{equi1} we have
$\mbox{div\,} {\bf j}=0$, then it is possible to work out the poloidal-toroidal decomposition for the current  ${\bf j}={\bf j}_{\perp}+{\bf j}_{\varphi}$ as well. In particular there exists a function $f$ such that $\mu_0{\bf j}_{\perp}=\frac{1}{r}\grad f \times {\bf e}_{\varphi}$ and since from \rife{equi1} we have ${\bf j}\cdot \grad p=0$, then $\grad f \times \grad p=0$.
Therefore $f$ is also constant on a magnetic surface and we can assume $f=f(\psi)$.\\

At this point we deduce from \rife{equi2} that $B_{\varphi}=\frac{f}{r}{\bf e}_{\varphi}$ and
\beq\label{equi4}
-\Delta^{*}\psi =  \mu_0j_{\varphi}\mbox{ in }\om_{+},
\eeq
where $\Delta^{*}=\nabla \left(\frac{1}{r} \nabla\right)=
\frac{\pa }{\pa r}\left(\frac{1}{r}\frac{\pa }{\pa r}\right)+\frac{1}{r}\frac{\pa^2}{\pa z^2}$.
Next observe that,
$$
({\bf j}_{\perp}+j_{\varphi}{\bf e}_{\varphi})\times
({\bf B}_{\perp}+B_{\varphi}{\bf e}_{\varphi})=
-\frac{1}{\mu_0 r}B_{\varphi}\grad f + \frac{1}{ r}j_{\varphi}\grad \psi=
-\frac{1}{\mu_0 r^2}f\grad f- \frac{1}{\mu_0 r}\Delta^{*}\psi \grad \psi.
$$
By using $\grad f =f^{'}\grad \psi$, $\grad p =p^{'}\grad \psi$ and
substituting in \rife{equi1} we conclude that,
$$
\left(- \Delta^{*}\psi-\frac{1}{r}ff^{'}-\mu_0rp^{'}\right)\grad \psi=0,
$$
which we solve setting
\beq\label{GS1}
- \Delta^{*}\psi=\frac{1}{ r}ff^{'}+\mu_0rp^{'}  \mbox{ in }\om.
\eeq

However recall from \rife{equi4} that the right hand side of \rife{GS1} is $j_{\varphi}$ in $\om_+$,
whence more exactly we have
\beq\label{GS2}
- \Delta^{*}\psi=\frac{1}{r}ff^{'}+\mu_0 rp^{'}\; \mbox{ in }\om_{+},
\quad - \Delta^{*}\psi=0\; \mbox{ in }\om_{-},
\eeq
which is the Grad-Shafranov equilibrium equation.\\

Concerning the boundary conditions let ${\nu}$ and ${\tau}$ denote unit outer normal and tangent
vector respectively to $\pa \om$ or $\pa \om_{+}$. Then we impose the natural boundary
conditions,
$$
{\bf B}_{\perp}\cdot{\bf \nu}=0 \mbox{ on } \pa \om
$$
$$
{\bf B}_{\perp}\cdot{\bf \nu}=0 \mbox{ and } {\bf B}_{\perp}\cdot{\bf \tau} \mbox{ is continuous } \mbox{ on } \pa \om_{+}
$$
and then, because of \rife{bdy1}, \rife{bdy2}, we also have,
$$
{\bf B}_{\perp}\cdot{\bf \nu}=\frac{1}{r}\frac{\pa \psi }{\pa \tau},\quad
{\bf B}_{\perp}\cdot{\bf \tau}=\frac{1}{r}\frac{\pa \psi }{\pa \nu}.
$$
Therefore $\frac{\pa \psi }{\pa \tau}=0$ locally on $\pa\om_+$ and $\pa\om$, whence $\psi$ is single valued in $\om$, constant on $\pa \om_{+}$ and $\pa \om$ and we fix
\beq\label{bdy3}
\psi=\gamma \mbox{ on } \pa\om.
\eeq
Also $\frac{\pa \psi }{\pa \nu}$ must be continuous on $\om_{+}$ and then, since $\Delta^{*}=\nabla \left(\frac{1}{r} \nabla\right)$, then the total
current $I$ of the plasma satisfies,
\beq\label{bdy4}
I=\mu_0\int\limits_{\om_+}j_{\varphi}=-\int\limits_{\om_+}\Delta^*\psi=-\int\limits_{\pa \om_+}\frac{1}{r}\frac{\pa \psi }{\pa \nu}=\int\limits_{\pa \om_-}\frac{1}{r}\frac{\pa \psi }{\pa \nu}=-\int\limits_{\pa \om}\frac{1}{r}\frac{\pa \psi }{\pa \nu},
\eeq
where we used that $\Delta^*\psi=0$ on $\om_{-}$.\\

At last, since $\om$ is a bounded, smooth simply connected domain whose closure is contained in the halfplane $\{x>0\}$, then the operator $\Delta^*$ is smooth and uniformly elliptic in $\om$ and setting
$(x_1,x_2)=(r,z)$ and
$g(x_1,x_2,\psi)=\frac{1}{ x_1}f(\psi) f^{'}(\psi)+\mu_0 x_1 p^{'}(\psi)$ we have that \rife{GS2} is the same as
\beq\label{GS4}
- \Delta^* \psi=[g(x_1,x_2,\psi)]_{+}\mbox{ in }\om,
\eeq
where $[g]_+=\max\{g,0\}$, which is to be solved together with the boundary conditions \rife{bdy3},\rife{bdy4}.\\ \\

The exact forms of $p^{'},f,f^{'}$ and consequently of $g$ are unknowns of the problem and we refer to
\cite{blum} and the references quoted therein for a detailed historical account of the so called "reconstruction problem" of these structure functions. As a consequence the analytical results used in the physical literature are often dictated by various assumptions and approximations depending by the specific problem at hand, see \cite{CFrei, Dubin, Fre, Porc} and the references quoted therein.\\

This was the initial motivation to pursue the rigorous analysis of \rife{GS4} with boundary conditions \rife{bdy3}, \rife{bdy4} under very general assumptions both on the uniformly elliptic operator and on $g$ (\cite{AmbM,BMar1,BeBr,Te}). However, while general existence results are available (\cite{BeBr,Te}),
the understanding of the more physically relevant global structure of the solutions set and the related bifurcation diagram is in general a challenging open problem. In particular, uniqueness is known in general only for $g$ uniformly Lipschitz in $\psi$ (\cite{BeBr,Te2}), whence of at most linear growth.
This is our first motivation to tackle a detailed study of the bifurcation diagram for \fbi\,, that is, in the more involved case where $g$ is of superlinear growth, while, as a standard simplifying assumption (the large aspect ratio approximation \cite{Fre,Taylor}), the operator $\Delta^*$ is replaced by the standard Laplace operator $\Delta=\frac{\pa^2}{\pa x_1^2}+\frac{\pa^2}{\pa x_2^2}$. Of course, we think at this as a first step toward the analysis of the more general problem \rife{bdy3}, \rife{bdy4}, \rife{GS4}. Besides, \fbi\, naturally arises
by a suitably defined variational principle (\cite{BeBr}, see also Appendix A below)
in the same spirit of the well known mean field theories classically adopted in the statistical mechanics description of bounded plasmas and vortex systems (\cite{clmp2, MoJo, MoTV, SmithO, Taylor}). These entropy maximization principles are concerned with two-dimensional plasmas described by a problem of the form \fbi\, (with
exponential nonlinearity) whence suitable to describe the above mentioned
large aspect ratio limit.\\

For the model problem \fbi\, a uniqueness result has been recently obtained by the first and third author
up to a certain threshold (\cite{BJ1}, see Theorem B below).
In this paper we make a further progress in that direction and provide an alternative about the global bifurcation diagram, see Theorem \ref{alt}. It is well known that Tokamaks with vertically elongated and D-shaped cross sections have better performances than the standard ones with circular shapes, see \cite{hoff, Liue} and references therein. Therefore, it seems also interesting that the alternative about the bifurcation diagram presented here holds
for "most domains" in a suitably defined sense, see Definition \ref{D1}.\\
Another result of physical relevance is (see Theorem \ref{thlambda})
that for high enough values of the current $I$ it happens that $\gamma<0$, showing that indeed the plasma does not touch the boundary of the vessel in this range.

\bigskip
\bigskip

We will focus here on positive solutions of \fbi. To explain our goals we first move to the well known dual formulation (\cite{BeBr,Te2}) of \fbi\, which, for positive solutions, consists of the constrained problem,
$$
\graf{-\Delta \psi =(\al+{\lm}\psi)^p\quad \mbox{in}\;\;\om\\ \\
\bigintss\limits_{\om} {\dsp \left(\al+{\lm}\psi\right)^p}=1\\ \\
\psi>0 \quad \mbox{in}\;\;\om, \quad \psi=0 \quad \mbox{on}\;\;\pa\om \\ \\
\al\geq0
}\qquad \prl
$$
for the unknowns $\al\in\R$ and $\psi \in C^{2,r}_{0,+}(\ov{\om}\,)$.
Here, $\lm\geq 0$ and $p\in [1,p_{N})$ are fixed and for $r\in (0,1)$ we set
$$
C^{2,r}_0(\ov{\om}\,)=\{\psi \in C^{2,r}(\ov{\om}\,)\,:\, \psi=0\mbox{ on }\pa \om\},\;
C^{2,r}_{0,+}(\ov{\om}\,)=\{\psi \in C^{2,r}_0(\ov{\om}\,)\,:\, \psi> 0\mbox{ in } \om\}.
$$

We remark that, as far as we assume $|\om|=1$, and since $\lm\geq 0$ by assumption, then if $(\all,\pl)$ solves $\prl$
we necessarily have,
$$
\all\leq 1,
$$
and the equality holds if and only if $\lm=0$. We will frequently use this fact without further comments.\\

We define positive solutions as follows.

\medskip

{\bf Definition.}
{\it We say that $(\gal,\val)$ is a positive solution
of \fbi\, if $\gal> 0$.}
{\it We say that $(\all,\pl)$ is a positive solution of {\rm $\prl$} if $\all>0$.}

\brm\label{remeq} {\it
Let $q$ be the conjugate index of $p$, that is
$$
\frac1p+\frac1q=1.
$$
For any fixed $\lm>0$ and $p>1$, $(\all,\pl)$ is a
solution of {\rm $\prl$} if and only if,
for $I=I_{\ssl}=\lm^{q}$, $(\gal,\val)=(\lm^{\frac{1}{p-1}}\all,\lm^{\frac{1}{p-1}}(\all+\lm\pl))$ is a
solution of \fbi. Therefore,
in particular, if $(\gal,\val)$ is a solution of \fbi\, with $\ga_{\sscp I}\geq 0$, then
$(\all,\pl)=(I^{-\frac{1}{p}}\gal,I^{-1}(\val -\gal))$ solves
{\rm $\prl$} and the identity $I^{-\frac{1}{p}}\val=\all+\lm\pl$ holds.}
\erm

\bigskip

It is well known (\cite{BeBr}) that for $p\in[1,p_{N})$ a solution of \fbi\, exists for any $I>0$.
In particular solutions of \fbi\, are found in \cite{BeBr,Te2} also as solutions of
suitably defined variational problems. For $p>1$ there exists a one to one correspondence (\cite{BeBr})
between variational solutions of \fbi\, with $\ga_{\sscp I}\geq 0$
and variational solutions of $\prl$, see Appendix \ref{appD}  for further details.
Actually, concerning variational solutions, as a consequence of the results in \cite{BMar1} (see also {\cite{BJ1}} and references therein), the following holds true.\\

{\bf Theorem A.} (\cite{BMar1})\\
{\it Let $p\in(1,p_{N})$ and $(\gal,\val)$ be a variational solution of \fbi.
There exists $I^{**}(\om,p)\in (0,+\ii)$ such that $\gal>0$ if and only if $I\in (0,I^{**}(\om,p))$
and $\gal=0$ if and only if $I=I^{**}(\om,p)$.\\
For $p\in(1,p_{N})$ define $\lm^{**}(\om,p)=(I^{**}(\om,p))^{\frac1q}$ and let $(\all,\pl)$ be a
variational solution of {\rm $\prl$}.
Then $\all>0$ if and only if $\lm\in [0,\lm^{**}(\om,p))$ and $\all=0$ if and only if $\lm=\lm^{**}(\om,p)$.
}

\

In particular Theorem A shows that for any $\lm<\lm^{**}(\om,p)$ there exists at least
one (variational) solution of $\prl$. The first and third author have recently proved
a uniqueness and monotonicity result (\cite{BJ1}).
For fixed $t\geq 1$ we denote
\beq\label{Sob}
\Lambda(\om,t)=\inf\limits_{w\in H^1_0(\om), w\equiv \!\!\!\!/ \;0}\dfrac{\ino |\nabla w|^2}{\left(\ino |w|^{t}\right)^{\frac2t}}\,,
\eeq
which provides the best constant in the Sobolev embedding $\|w\|_p\leq S_p(\om)\|\nabla w\|_2$, $S_p(\om)=\Lambda^{-2}(\om,p)$,
$p\in[1,2p_{N})$.
For $(\all,\pl)$ a solution of $\prl$ we define the energy,
$$
\el:=\frac12 \ino\rl \pl\equiv \frac12 \ino |\nabla \pl|^2,
$$
and
$$
\lm^*(\om,p)=\sup\Bigr\{\lm>0\,:\,\al_{\mu}>0\,\mbox{\rm \,for any solution of }
{\textbf{(}\mathbf P\textbf{)}_{\mathbf \mu}},\,\forall\,\mu<\lm\Bigr\},
$$
{where $\rl$ is defined later in \eqref{eq:rl}}. We denote by $\mathbb{D}_N$ the ball of unit volume. {For $|\om|=1$,} let $\mathcal{G}(\om)$ denote the set of solutions $(\all,\pl)$ of $\prl$ in
$[0,\frac1p \Lambda(\om,2p))$. For the sake of clarity, we point out that if a map $\mathcal{M}$
from an interval $[a,b]\subset \R$ to a Banach space $X$ is said to be real analytic,
then it is understood that $\mathcal{M}$ can be extended in an open neighborhood
of $a$ and $b$ where it admits locally a power series expansion, totally convergent in the $X$-norm.\\

{Concerning the unique branch of positive solutions, the first author and third author recently proved the following result:}

\smallskip

{\bf Theorem B.} {(\cite{BJ1})}
{\it
Let $\om\subset \R^N$ be a bounded domain of class $C^3$, $|\om|=1$ and $p\in [1,p_{\sscp N})$. Then $\lm^*(\om,p)\geq \frac1p \Lambda(\om,2p)$ and the equality holds if and only if $p=1$. Moreover, we have:
\begin{itemize}
\item[1.] \emph{(Uniqueness):} for any $\lm \in [0,\frac1p \Lambda(\om,2p))$
there exists a unique solution $(\all,\pl)$ of {\rm $\prl$}.

\item[2.] \emph{(Monotonicity):} $\mathcal{G}(\om)$ is a real analytic simple curve of positive solutions $[0,\frac1p\Lambda(\om,2p))\ni \lm\mapsto (\all,\pl)$ such that, for any $\lm\in [0,{\textstyle \frac1p}\Lambda(\om,2p))$,
$$
\frac{d \all}{d\lm}<0\quad\mbox{ and }\quad \frac{d \el}{d\lm}>0,
$$
and
$$
\all=1+\mbox{\rm O}(\lm),\;\;\pl=\psi_{\sscp 0}+\mbox{\rm O}(\lm),
\;\; \el=E_0(\om)+\mbox{\rm O}(\lm), \quad \mbox{ as }\lm\to 0^+,
$$
where,
$$
E_{0}(\om)=\frac12\ino \ino G_{\om}(x,y)\,dxdy\leq
E_{0}(\mathbb{D}_{\sscp N})=
{\textstyle\dfrac{|B_1|^{\scp -\frac 2 N}}{4(N+2)}}.
$$
In particular $\mathcal{G}(\om)$ can be extended continuously on $\lm \in [0,\frac1p \Lambda(\om,2p)]$ with
$(\all,\pl)\to (\ov{\al},\ov{\psi})$ as $\lm\to \frac1p\Lambda(\om,2p)^-$ and $\ov{\al}=0$ if and only if $p=1$.
\end{itemize}
}

\

Observe that from Theorems A and B we have $\lm^{**}(\om,p)\geq \lm^*(\om,p)>\frac1p\Lambda(\om,2p)$, for $p>1$.
On the other hand, the set of positive variational solutions of $\prl$ is
not empty for any $\lm\in (\frac1p\Lambda(\om,2p), \lm^{**}(\om,p))$, see \cite{BMar,BeBr}.
Our main concern here is with respect to the continuation of the curve of solutions $\mathcal{G}(\om)$, under generic assumptions, beyond $\frac1p\Lambda(\om,2p)$,
enjoying uniqueness and monotonicity properties.\\

Our first result is about the existence of a free boundary
in the interior of $\om$ for solutions of \fbi\, with $p\in (1,p_{N})$ (the case $p=1$ is fully understood,
see \cite{BeBr,Pudam,Te2}).
The point here is that one would like to know whether or not, for a given $I$, $\ga_{\sscp_I}$ is negative,
which implies in particular that $\om_{-}:=\{x\in\om\,:\,v_{\sscp I}<0\}$ is not empty.
For $p>1$ the existence of a multiply connected free boundary has been proved in \cite{We} for $I$ large
and under some assumptions about the existence of non degenerate critical points of a suitably defined
{Kirchoff-Routh} type functional. Still for $I$ large, but only for $N=2$ and for domains with non trivial topology,
a similar result has been obtained in \cite{Liu}.
Other sufficient conditions for the existence of solutions with $\ga<0$ has been found in \cite{AmbM},
which however assume the nonlinearity $v_+^p$ to be replaced by $g_+(x,v)$ satisfying $g(x,t)\geq ct$,
for some $c>0$, which therefore does not fit our problem.
As mentioned above, for variational solutions we have $\ga<0$ if and only if $I>I^{**}(\om,p)$ (\cite{BMar1}).

We prove here that in fact this is always true, in the sense that
for any $\om$ and $p\in (1,p_N)$ we have that for any $I$
large enough it holds $\ga< 0$. More {exactly} we show that there are no
solutions of $\prl$ with $\all\geq 0$ for $\lm$ large.
\bte\label{thlambda} Let $\om\subset \R^N$ be a bounded domain of class $C^3$
and $p\in(1,p_N)$.

\smallskip

$(a)$ There exists
$\ov{\lm}=\ov{\lm}(\om,p)>0$ depending only on $p, N$ and $\om$ such that if $(\lm,\pl)$ is a
solution of {\rm$\prl$}, then $\lm\leq \ov{\lm}$.\\
$(b)$ Let $\ov{I}=\ov{I}(\om,p)=\ov{\lm}^{\frac{p}{p-1}}(\om,p)$, then $\gamma_{\sscp I}< 0$ for any solution
$(\gal,\val)$ of \fbi\, with $I>\ov{I}$.
\ete

\medskip

Observe that the existence of a free boundary is obtained with no assumptions on the domain.

\bigskip

We come now to the continuation of the curve of solutions $\mathcal{G}(\om)$. It has been shown in \cite{BJ1} that if
$\om=\mathbb{D}_2$, then $\mathcal{G}(\mathbb{D}_2)$ can be continued on the maximal
interval $\lm\in [0,\lm^*(\mathbb{D}_2,p)]$ where $\all$ is monotonic decreasing and $\all\to 0^+$ as
$\lm\to \lm^*(\mathbb{D}_2,p)^{-}$ while $\el$ is monotonic increasing and
$\el\to E_*(\mathbb{D}_2,p):=\frac{p+1}{16\pi}$
as $\lm\to \lm_*(\mathbb{D}_2,p)^{-}$. Interestingly enough, for $N=2$ and $|\om|=1$ one always has $E_{\lm}(\om,p)\leq E_*(\mathbb{D}_2,p)$, as recently shown in \cite{BJ3}. The upper bound is optimal and the equality holds if and only if,
up to a translation, $\om=\mathbb{D}_2$. The value of $\lm^*(\mathbb{D}_2,p)$ has also been
evaluated explicitly in \cite{BJ3}, please see $(j)$ in Theorem \ref{lmconvex} below.
However, it is not trivial to come up with the maximal branch $\mathcal{G}(\mathbb{D}_2)$  as one needs
to rule out the existence of bifurcation points as well as the monotonicity of $\all,\el$ for
$\lm>\frac1p\Lambda(\om,2p)$. In particular the uniqueness of solutions
on $\mathbb{D}_2$ (\cite{BSp}) was used in a crucial way in \cite{BJ1}, which also implies that
$\ov{\lm}(\mathbb{D}_2,p)$, as defined by Theorem \ref{thlambda}, satisfies
$\ov{\lm}(\mathbb{D}_2,p)=\lm^*(\mathbb{D}_2,p)$.
In general the situation is much more difficult and we
come up with an alternative, in a suitably defined generic sense.\\

For $\om_0$ a bounded domain of class $C^4$,
we denote by $\mbox{\rm Diff}^{\,4}(\om_0)$ the set of domains $\om$ such that $\om=h(\om_0)$ for some diffeomorphism
$h:\ov{\om_0}\to\ov{\om}$ of class $C^4$. We recall that a subset of a metric space is said to be:\\
- nowhere dense, if its closure has empty interior; \\
- meager (or of first Baire category), if it is the union of countably many nowhere dense sets.

\smallskip

\bdf\label{D1}
{\it We say that a property holds for most bounded domains $\om$ of class $C^4$, $|\om|=1$, if, given any bounded domain $\om_0$ of class $C^4$, it holds on $\om=|h(\om_0)|^{-\frac1N}h(\om_0)$ with $h\in \mbox{\rm Diff}^{\,4}(\om_0)\setminus\mathcal{F}$ for some meager set $\mathcal{F}\subset \mbox{\rm Diff}^{\,4}(\om_0)$.}
\edf

\

Then we have,
\bte\label{alt}
Let $p\in(1,p_N)$. For most bounded domains $\om\subset \R^N$ of class $C^4$, $|\om|=1$, see Definition \ref{D1}, it holds, either:

\medskip

$(i)$ there exists $\lm_\ii(\om,p)\in (\frac1p\Lambda(\om,2p),\ov{\lm}(\om,p)]$ such that
$\mathcal{G}(\om)$ can be continued to a real analytic simple curve of positive solutions
defined in $[0, \lm_\ii(\om,p))$ such that,
$$
\all\mbox{ is strictly decreasing and } \el\mbox{ is strictly increasing in }
(0,{\lm}_\ii(\om,p)),
$$
$$
\all\searrow 0^+,\quad \el\nearrow {E}_{\ii}(\om,p),
\mbox{ as }\lm \nearrow \lm_\ii(\om,p),
$$
and for any sequence $\lm_n\to \lm^{-}_\ii(\om,p)$ there exists a subsequence
$\{\lm_{\sscp k_n}\}\subseteq \lm_n$ such that $\psi_{\sscp \lm_{\sscp k_n}}\to \psi_{\ii}$
in $C_0^2(\ov{\om}\,)$, where $\psi_{\ii}$ solves {\rm $\prl$} with $\lm=\lm_\ii(\om,p)$
and $\al=0$ and ${E}_{\ii}(\om,p)=\left. \el\right|_{\lm={\lm}_\ii(\om,p)}$, or

\

$(ii)$ there exists $\lm_\ii(\om,p)\in (\frac1p\Lambda(\om,2p),\ov{\lm}(\om,p)]$ such that
$\mathcal{G}(\om)$ can be continued to a continuous simple curve without bifurcation points,
$$
\mathbb{G}_\ii=\Bigr\{[0,s_{\ii})\ni s\mapsto (\lm(s),\al(s),\psi(s))\Bigr\}
$$
which has locally also the structure of a $1$-dimensional real analytic manifold,
such that for any $s$, $(\al(s),\psi(s))$ is a positive solution of
{\rm $\left.\prl\right|_{\lm=\lm(s)}$} and, as $s\to s_{\ii}$, we have $\al(s)\to 0^+$,
$\lm(s)\to \lm_\ii(\om,p)$ and
for any sequence $s_n\to s_\ii$, there exists a subsequence $\{t_n\}\subseteq \{s_n\}$,
such that $\psi(t_n)\to \psi_{\ii}$ in $C_0^2(\ov{\om}\,)$, which solves {\rm $\prl$} with $\lm=\lm_\ii$ and $\al=0$.
\ete

\bigskip

\brm {\it Interestingly enough, not only in a generic sense either $(i)$ or $(ii)$ hold
but also the $\lm$-limit set of the curve of solutions on the section $\al=0$ is a singleton.
Besides, the proof shows in particular that it is well defined,
$$
{\lm}_1={\lm}_1(\om):=\sup\left\{\mu>\frac1p\Lambda(\om,2p)\,:\,\all>0\mbox{ and }\sg_1(\all,\pl)>0,\;
\forall\,(\all,\pl)\in \mathcal{G}_{\mu},\,\forall \lm<\mu\right\},
$$
where $\mathcal{G}_{\mu}$ is the continuation of $\mathcal{G}(\om)$, $\sg_1(\all,\pl)$ is a suitably
defined first eigenvalue (see {\rm \rife{4.1}} below) and in particular that if $\al_{\sscp \lm_1(\om)}=0$, then $(i)$ holds with $\lm_\ii(\om,p)=\lm_1(\om)$.
It is a challenging open problem to understand under which conditions on $\om$ this property holds.}
\erm

\medskip

Let us point out that by Theorem 1.4 in \cite{BJ1} we know that $(i)$ of Theorem \ref{alt} holds for $\om=\mathbb{D}_2$. It turns out that under a natural connectivity assumption about the set of variational solutions, $(i)$ of Theorem \ref{alt} always holds. We refer to Section 5 for full details.

\medskip

Moreover, for symmetric and convex domains the above generic property is fulfilled and then either $(i)$ or $(ii)$ of Theorem \ref{alt} hold. Indeed, let $\om\subset \R^2$ be symmetric and convex with respect to the coordinate
directions or $\om=B_{R}\subset \R^{N}$, $N\geq 3$. By some uniqueness results in \cite{DGP} we know that there exists a unique value of $\lm$, which we denote by
$\lm_*(\om,p)$, such that there exists a solution $\psi_*$ of {\rm $\prl$} with $\al=0$, see Theorem C at the end of Section 4.
Then we have,
\bte\label{lmconvex}$\,$\\
If $\om\subset \R^2$, $|\om|=1$, is symmetric and convex with respect to the coordinate
directions or if $\om=\mathbb{D}_N$, $N\geq 3$,
then either $(i)$ or $(ii)$ of Theorem \ref{alt} hold and $\lm_\ii(\om,p)=\lm_*(\om,p)$, $\psi_\ii=\psi_*$.\\
Moreover, it holds:

\medskip

$(j)$ If $\om\subset \R^2$, $|\om|=1$, is as in Theorem \ref{alt}, then
$$\lm^{2p}_\ii(\om,p)\geq \left(\frac{8\pi}{p+1}\right)^{p-1}\Lambda^{p+1}(\om,p+1)$$
and the equality
holds if and only if up to a translation, $\om=\mathbb{D}_2$.

\medskip

$(jj)$ If $\om\subset \R^2$, $|\om|=1$, is symmetric and convex with respect to the coordinate
directions, then
\beq\label{030221.1}
\left(\frac{8\pi}{p+1}\right)^{p-1}\Lambda^{p+1}(\om,p+1)\leq \lm^{2p}_*(\om,p)\leq
\frac{\left(\frac{p+1}{8\pi}\right)^{p+1}}{(2 E_*(\om,p))^{2p}}
\Lambda^{p+1}(\om,p+1)
\eeq
where $E_*(\om,p)$ is the energy
of the solution $\psi_*$. Moreover, the l.h.s. equality holds if and only if, up to a translation, $\om=\mathbb{D}_2$.
If $\om=\mathbb{D}_2$ the r.h.s. inequality is an equality.
\ete

\medskip

The last estimate is a refinement of Theorem 1.2
in \cite{BJ3}, which states that if $\om\subset \R^2$, $p>1$ and $(\lm,\pl)$ solves {\rm $\prl$}
with $\all=0$, then
$\lm^{2p}\geq \left(\frac{8\pi}{p+1}\right)^{p-1}\Lambda^{p+1}(\om,p+1)$,
where the equality holds if and only if, up to a translation, $\om=\mathbb{D}_2$.

\

Concerning Theorem \ref{alt}, we handle the problem of possible
bifurcation points by a generalization (see Proposition 2.7 in \cite{BJ1}) of the
Crandall-Rabinowitz (\cite{CrRab}) continuation Theorem to the constrained problem $\prl$.
A subtle point arises as we need a suitable transversality condition to be verified, which allows one to
prove that critical values of $\lm$ are not bifurcation points.
The fact that the tranversality condition is generic is not at all trivial
and is proved by a careful adaptation to the constrained problem $\prl$ of the
refined and simplified transversality theory with respect to domains variations
developed in \cite{H}.

\

Finally, we present an application of our approach for the problem
$$
\graf{-\Delta u=\mu (1+u)^p \quad \mbox{in}\;\;\om\\ u=0 \quad \mbox{on}\;\;\pa \om}
$$
which, interestingly enough, yields a sort of refinement of a result in \cite{GM} about its Rabinowitz unbounded continuum of solutions, see Section 6 for full details.

\bigskip
\bigskip

{This paper is organized as follows: Section 2 includes the functional set-up and preliminaries of the bifurcation analysis which will be essentially used later in the proof of Theorem \ref{alt}; Section~3 is devoted to the proof of Theorem \ref{thlambda} concerning the existence of a free boundary; in Section~4 we derive Theorems \ref{alt} and \ref{lmconvex} about the continuation of the branch of solutions, where a general transversality theorem comes into play; in Section 5, we refine the result of Theorem~\ref{alt} by imposing a natural connectedness assumption; a discussion on the application of this global bifurcation framework introduced here to a nonlinear eigenvalue problem is presented in Section~6.}

\bigskip
\bigskip

\section{\bf Preliminary results}\label{sec2}
In this section we report some results in \cite{BJ1} about the
spectral and bifurcation analysis for solutions of $\prl$ with $\lm\geq 0$ and
$p\geq 1$.\\
We write $\rla=\rla(\psi)=(\al+\lm\psi)^p$ and
{\begin{equation}\label{eq:rl}
\rl=(\all+\lm\pl)^p,
\end{equation}} and then in particular,
$$
(\rla(\psi))^{\frac1q}=(\al+\lm\psi)^{p-1}\mbox{ and }\rlq=(\all+\lm\pl)^{p-1},\mbox{ respectively},
$$
where it is understood that $(\rla(\psi))^{\frac1q}\equiv 1 \equiv \rlq$ for $p=1$.
Whenever $(\all,\pl)$ is a solution of $\prl$ and $\{\eta,\varphi\}\subset L^2(\om)$ we set,
$$
<\eta>_{\ssl}=\dfrac{\ino \rlq  \eta}{\ino \rlq},\quad \mbox{and }\quad  [\eta]_{\ssl}=\eta \,-<\eta>_{\ssl},
$$
and then we define,
$$
<\eta,\varphi>_{\ssl}:=\dfrac{\ino \rlq  \eta\varphi}{\ino \rlq}\quad \mbox{and}\quad
\|\varphi\|_{\ssl}^2:=<\varphi,\varphi>_{\ssl}=<\varphi^2>_{\ssl}=\dfrac{\ino \rlq  \varphi^2}{\ino \rlq}.
$$
Clearly for non negative solutions $\rl$ is strictly positive in ${\om}$.
Therefore it is easy to see that $<\cdot,\cdot>_{\ssl}$ defines a scalar product on $L^2(\om)$ whose norm is  $\|\cdot\|_{\ssl}$.
We will also adopt when needed the useful shorthand notation,
\beq\label{emel}
\ml=\ino \rlq.
\eeq

\brm\label{rem1}{\it
We will use the fact that,{\rm
$$
{\ino \rlq  [\eta]_{\ssl}^2}=\ino \rlq  (\eta\,-<\eta>_{\ssl})^2\geq 0,
$$
}where the equality holds if and only if $\eta$ is constant, whence in particular, in case $\eta\in H^{1}_0(\om)$,
if and only if $\eta$ vanishes identically. Also, since obviously $<[\eta]_{\ssl}>_{\ssl}=0$, then,
$$
<[\varphi]_{\ssl} \eta>_{\ssl}=<\varphi, [\eta]_{\ssl}>_{\ssl}= <[\varphi]_{\ssl}[\eta]_{\ssl}>_{\ssl},\;\fo\,\{\eta,\varphi\}\subset L^2(\om).
$$
We will use these properties time to time when needed without further comments.}
\erm

\bigskip
In the sequel we aim to describe branches of solutions of $\prl$ around a positive solution, i.e.
with $\all>0$. To this end we will use the following:
\brm\label{Analytic}
{\it It is not difficult to construct an open subset $A_{\sscp \om}$ of the
Banach space of triples $(\lm,\al,\psi)\in \R\times \R\times  C^{2,r}_0(\ov{\om}\,)$
such that, on $A_{\sscp \om}$, the density {\rm $\rla=\rla(\psi)=(\al+\lm\psi)^p$} is well defined and
$$
\all+\lm\pl\geq \frac{\all}{2} \quad \mbox{ in }\quad \ov{\om},
$$
in a sufficiently small open neighborhood in $A_{\sscp \om}$ of any triple of the form $(\lm,\all,\pl)$ whenever $(\all,\pl)$ is a positive solution of {\rm $\prl$}. As a consequence, relying on
known techniques about real analytic functions on Banach spaces {\rm (}\cite{but}{\rm )},
one can see that the map $\Phi(\lm,\al,\psi)$ below is jointly real analytic in an open neighborhood of $A_{\sscp \om}$
around any such triple $(\lm,\all,\pl)$.}
\erm

At this point we can introduce the maps,
\beq\label{eF}
F: A_{\sscp \om} \to  C^{r}(\ov{\om}\,),\quad F(\lm,\al,\psi):=-\Delta \psi -\rla(\psi),
\eeq

\beq\label{ePhi}
\Phi:A_{\sscp \om} \to  \R\times C^{r}(\ov{\om}\,),\quad
\Phi(\lm,\al,\psi)=\left(\begin{array}{cl}F(\lm,\al,\psi)\\ \\-1+\ino \rla\end{array}\right),
\eeq
and, for a fixed $(\lm,\al,\psi)\in A_{\sscp \om}$,
its differential with respect to $(\al,\psi)$, that is the linear operator,
$$
D_{\al,\psi}\Phi(\lm,\al,\psi):\R \times C^{2,r}_0(\ov{\om}\,) \to  \R\times C^{r}(\ov{\om}\,),
$$
which acts as follows,
$$
D_{\al,\psi}\Phi(\lm,\al,\psi)[s,\phi]=
\left(\begin{array}{cr}
D_\psi F(\lm,\al,\psi)[\phi]+d_{\al}F(\lm,\al,\psi) [s]\\ \\
\ino \left(D_\psi \rla[\phi]+d_{\al}\rla [s]\right)\end{array}\right),
$$
where we have introduced the differential operators,
\beq\label{lin}
D_\psi F(\lm,\al,\psi)[\phi]=
-\Delta \phi -\lm p \rlqa \phi, \quad \phi \in C^{2,r}_0(\ov{\om}\,),
\eeq

\beq\label{lin1}
D_\psi \rla[\phi]=\lm p\rlqa \phi, \quad \phi \in C^{2,r}_0(\ov{\om}\,),
\eeq
and

\beq\label{lin2}
d_{\al}F(\lm,\al,\psi) [s]=-p \rlqa s, \quad s \in \R,
\eeq

\beq\label{lin3}
d_{\al}\rla [s]= p\rlqa s, \quad s \in \R.
\eeq

\bigskip

For fixed $\lm\geq 0$ and $p\geq 1$, the pair $(\all,\pl)$ solves {\rm$\prl$} if and only if  $\Phi(\lm,\all,\pl)=(0,0)$ and we define the linear operator,
\beq\label{eLl}
L_{\ssl}[\phi]=D_\psi F(\lm,\all,\pl)[\phi]=-\Delta \phi-\lm p \rlq [\phi]_{\ssl}.
\eeq

We say that $\sg=\sg(\all,\pl)\in\R$ is an eigenvalue of $L_{\ssl}$ if the equation,
\beq\label{lineq0}
-\Delta \phi-\lm p \rlq [\phi]_{\ssl}=\sg\rlq [\phi]_{\ssl},
\eeq
admits a non trivial weak solution $\phi\in H^1_0(\om)$. Let us define the Hilbert space,
\beq\label{y0}
Y_0:=\left\{ \varphi \in \{L^2(\om),<\cdot,\cdot >_{\ssl}\}\,:\,<\varphi>_{\ssl}=0\right\},
\eeq
and $T(f):=G[ \rlq f]$, for $f\in L^2(\om)$.
Since $T(Y_0)\subset W^{2,2}(\om)$, then the linear operator,
\beq\label{T0}
T_0:Y_0\to Y_0,\;T_0(\varphi)=G[\lm p \rlq\varphi]-<G[\lm p \rlq\varphi]>_{\ssl},
\eeq
is compact. By a straightforward evaluation we see that $T$ is also self-adjoint. As a consequence,
standard results concerning the spectral decomposition of self-adjoint, compact, linear operators on
Hilbert spaces show that $Y_0$ is the Hilbertian direct sum of the eigenfunctions of $T_0$, which can be represented as
$\varphi_k=[\phi_k]_{\ssl}$, $k\in\N=\{1,2,\cdots\}$,
$$
Y_0=\overline{\mbox{Span}\left\{[\phi_k]_{\ssl},\;k\in \N\right\}},
$$
for some $\phi_k\in H^1_0(\om)$, $k\in\N=\{1,2,\cdots\}$. In fact, the eigenfunction $\varphi_k$,
whose eigenvalue is $\mu_k\in \R\setminus\{0\}$, satisfies,
$$
\mu_k\varphi_k=\left(G[\lm p \rlq\varphi_k]-< G[\lm p \rlq\varphi_k]>_{\ssl}\right).
$$
In other words, by defining,
$$
\phi_k:=(\lm p+\sg_k) G[\rlq\varphi_k],
$$
it is easy to see that $\varphi_k$ is an eigenfunction of $T_0$ with eigenvalue $\mu_k=\frac{\lm p}{\lm p+\sg_k}\in \R\setminus\{0\}$ if and only if
$\phi_k\in H^1_0(\om)$ and weakly solves,
\beq\label{lineqm}
-\Delta \phi_k= (\lm p+\sg_k)  \rlq [\phi_k]_{\ssl}\quad\mbox{ in }\quad \om.
\eeq

In particular we will use the fact that $\varphi_k=[\phi_k]_{\ssl}$ and
\beq\label{lineq9m}
\phi_k=(\lm p+\sg_k) G[\rlq[\phi_k]_{\ssl}],\quad k\in\N=\{1,2,\cdots\}.
\eeq

At this point, standard arguments in the calculus of variations show that,
\beq\label{4.1}
\sg_1=\sg_1(\all,\pl)=\inf\limits_{\phi \in H^1_0(\om)\setminus \{0\}}
\dfrac{\ino |\nabla \phi|^2 - \lm p \ino \rlq [\phi]_{\ssl}^2 }{\ino \rlq [\phi]_{\ssl}^2}.
\eeq
The ratio in the right hand side of  \rife{4.1} is well defined because of
Remark \ref{rem1}. Higher eigenvalues are defined inductively via the variational problems,
\beq\label{4.1.k}
\sg_k(\all,\pl)=\inf\limits_{\phi \in H^1_0(\om)\setminus \{0\},<\phi,[\phi_m]_{\ssl}>_{\ssl}=0,\; m\in \{1,\ldots, k-1\}}
\dfrac{\ino |\nabla \phi|^2 - \lm p \ino \rlq [\phi]_{\ssl}^2 }{\ino \rlq [\phi]_{\ssl}^2}.
\eeq
The eigenvalues form a numerable non decreasing
sequence $\sg_1\leq \sg_2\leq....\leq \sg_k\leq...\,$. By the Fredholm alternative,
if $0\notin \{\sg_j\}_{j\in \N}$, then $I-T_0$ is an isomorphism of $Y_0$ onto itself.

\bigskip

Concerning $D_{\al,\psi}\Phi(\lm,\al,\psi)$ we have,

\bpr{{\rm (\cite{BJ1})}}\label{pr2.2} For any positive solution $(\all,\pl)$ of {\rm$\prl$} with $\lm\geq 0$, the kernel of
$D_{\al,\psi}\Phi(\lm,\all,\pl)$ is empty if and only if the equation, {\rm
\beq\label{contr}
-\Delta \phi-\lm p \rlq [\phi]_{\ssl}=0,\quad \phi\in C^{2,r}_0(\ov{\om}\,)
\eeq}
admits only the trivial solution, or equivalently, if and only if $0$ is not an eigenvalue of $L_{\ssl}$.
\epr

\bigskip

An easy to prove but relevant identity is satisfied by any eigenfunction, which we summarize in the following,
\ble{{\rm (\cite{BJ1})}}\label{lstr}
Let $(\all,\pl)$ be a solution of {\rm $\prl$} and let $\phi_k$ be any eigenfunction of an eigenvalue $\sg_k=\sg_k(\all,\pl)$.
Then the following identity holds,
\beq\label{2907.6n}
\frac{1}{\ml}<\phi_k>_{\ssl}\equiv(\all+\lm<\pl>_{\ssl})<\phi_k>_{\ssl}=(\lm(p-1)+\sg_k)<\pl [\phi_k]_{\ssl}>_{\ssl}.
\eeq
\ele

\bigskip

The implicit function theorem applies around a positive solution of $\prl$ in the following form:
\ble{{\rm (\cite{BJ1})}}\label{lem1.1} Let $(\al_{\sscp \lm_0},\psi_{\sscp \lm_0})$ be a positive solution of {\rm $\prl$} with $\lm=\lm_0\geq 0$.\\
If $0$ is not an eigenvalue of $L_{\sscp \lm_0}$, then:\\
$(i)$ $D_{\al,\psi}\Phi(\lm_0, \al_{\sscp \lm_0}, \psi_{\sscp \lm_0})$ is an isomorphism;\\
$(ii)$ There exists an open neighborhood $\mathcal{U}\subset A_{\sscp \om}$ of $(\lm_0,\al_{\sscp \lm_0},\psi_{\sscp \lm_0})$ such
that the set of solutions of
{\rm $\prl$} in $\mathcal{U}$ is a real analytic curve of positive solutions $J\ni\lm\mapsto (\all,\pl)\in B$, for
suitable neighborhoods $J$ of $\lm_0$ and $B$ of $(\al_{\sscp \lm_0},\psi_{\sscp \lm_0})$ in
$(0,+\ii)\times C^{2,r}_{0,+}(\ov{\om}\,)$.\\
$(iii)$ In particular if $|\om|=1$ and $(\al_{\sscp \lm_0},\psi_{\sscp \lm_0})=(\al_{\sscp 0},\psi_{\sscp 0})=(1,G[1])$, then
$(\all,\pl)=(1,\psi_{\sscp 0})+\mbox{\rm O}(\lm)$ as $\lm\to 0$.
\ele

\bigskip

Next we state a generalization of the bending result of \cite{CrRab} for solutions of $\prl$, which
is an improvement of Lemma \ref{lem1.1} in case of simple and vanishing eigenvalues satisfying a suitable transversality
condition.
\bpr{{\rm (\cite{BJ1})}}\label{pr3.1}
Let $(\all,\pl)$ be a positive solution of {\rm$\prl$} with $\lm>0$ and suppose that the  $k$-th
eigenvalue $\sg_k(\all,\pl)=0$ is simple, that is, it admits only one eigenfunction, $\phi_k\in C^{2,r}_0(\ov{\om}\,)$.
If $<\phi_k>_{\ssl}\neq 0$,
then there exists $\eps>0$, an open neighborhood $\mathcal{U}$ of $(\lm,\all,\pl)$ in $A_{\sscp \om}$ and
a real analytic curve $(-\eps,\eps)\ni s \mapsto (\lm(s), \al(s),\psi(s))$ such that
$(\lm(0), \al(0),\psi(0))=(\lm,\all,\pl)$ and the set of solutions of {\rm $\prl$} in $\mathcal{U}$ has the form $(\lm(s),\al(s),\psi(s))$,
where $(\al(s),\psi(s))$ is a solution of {\rm$\prl$} for $\lm=\lm(s)$ for any $s\in (-\eps,\eps)$, with
$\psi(s)=\pl+s\phi_k+\vxi(s)$, and
\beq\label{2907.0}
<[\phi_k]_{\sscp \lm(s)},\vxi(s)>_{\sscp \lm(s)}=0,\quad s\in (-\eps,\eps).
\eeq
Moreover it holds,
\beq\label{2907.1}
\vxi(0)\equiv 0\equiv \vxi^{'}(0),\quad \al^{'}(0)=-\lm <\phi>_{\ssl},\quad\lm^{'}(0)=0,\quad \psi^{'}(0)= \phi_k,
\eeq
and either $\lm(s)=\lm$ is constant in $(-\eps,\eps)$ or
$\lm^{'}(s)\neq 0$, $\sg_k(s)\neq 0$ in $(-\eps,\eps)\setminus\{0\}$, $\sg_k(s)$ is simple in $(-\eps,\eps)$ and

\beq\label{2907.5}
<[\phi_k]_{\ssl},\pl>_{\ssl}\neq 0\mbox{ and $<[\phi_k]_{\ssl},\pl>_{\ssl}$  has the same sign as }<\phi_k>_{\ssl},
\eeq

\beq\label{2907.10}
\dfrac{\sg_k(s)}{\lm^{'}(s)}=\dfrac{p<[\phi_k]_{\ssl},\psi_{\ssl}>_{\ssl}+\mbox{\rm o}(1)}
{<[\phi_k]_{\ssl}^2>_{\ssl}+\mbox{\rm o}(1)},\mbox{ as }s\to 0.
\eeq
\epr

\bigskip

The next result is about the transversality condition.
\bte\label{spectral1}{{\rm (\cite{BJ1})}} Let $(\all,\pl)$ be a positive solution of {\rm$\prl$}. Suppose that any eigenfunction $\phi_k$ of a fixed
vanishing eigenvalue $\sg_k=\sg_k(\all,\pl)=0$ satisfies $<\phi_k>_{\ssl}\neq 0$. Then $\sg_k$ is simple, that is,
it admits only one eigenfunction.\\
Let either $\om\subset \R^2$ be symmetric and convex with respect to the coordinate directions
$x_i$, $i=1,2$ or $\om=B_{R}\subset \R^{N}$, $N\geq 3$,
and let $(\all,\pl)$ be a positive solution of {\rm$\prl$} with $\lm>0$.
Suppose that $\sg_k=\sg_k(\all,\pl)=0$ and let $\phi_k$ be any
corresponding eigenfunction. Then:\\
$(i)$ $<\phi_k>_{\ssl}\neq 0$;\\
$(ii)$ $\sg_k(\all,\pl)$ is simple, that is, it admits at most one eigenfunction.
\ete

\brm {\it
Actually the proof of Theorem \ref{spectral1} in \cite{BJ1} covers only the case where
$\om\subset \R^2$ is symmetric and convex with respect to the coordinate directions.
However the proof is just an application of Theorem 3.1 in \cite{DGP} and it is
straightforward to check that if
$\om=B_{R}\subset \R^{N}$, $N\geq 3$ then the same argument works exactly in the same way just
by using Remark 3.1 in \cite{DGP}.}
\erm

\bigskip

The following proposition states a monotonicity property of $\all,\el$.

\bpr{{\rm (\cite{BJ1})}}\label{pr-enrg} Let $(\al_{\ssl_0},\psi_{\ssl_0})$ be a positive solution of {\rm $\prl$} with $\lm_0\geq 0$
and suppose that $0$ is not an eigenvalue of $L_{\ssl_0}$.
Then, locally near $\lm_0$, the map $\lm\mapsto (\all,\pl)$ is a real analytic simple curve of positive solutions
and if $\sg_1=\sg_1(\all,\pl)> 0$, then
$$
\frac{d\all}{d\lm}<0,\;\frac{d \el}{d\lm}>0.
$$
\epr

\bigskip
The following Lemma ensures that the first eigenvalue $\sg_1(\all,\pl)$ of a positive
variational solution must be non-negative.

\ble{{\rm (\cite{BJ1})}}\label{prmin}
Let $(\all,\pl)$ be a positive variational solution of {\rm $\prl$}. Then $\sg_1(\all,\pl)\geq 0$.
\ele

\bigskip

We will also need the following uniform estimate which is well known {{\rm (\cite{BeBr,BJ1})}}.
\ble\label{lemE1} Let $p\in[1,p_{N})$. For any $\ov{\lm}>0$ there
exists a positive constant $C_1=C_1(r,\om,\ov{\lm},p,N)$ depending only on $\om$, $\ov{\lm}$, $p$, $N$
and $r\in [0,1)$ such that $\|\pl\|_{C^{2,r}_0(\ov{\om})}\leq C_1$ for any solution
$(\all,\pl)$ of {\rm $\prl$} with $\lm\in [0,\ov{\lm}\,]$.
\ele

\bigskip
\bigskip

\section{\bf A uniform upper bound for $\lm$}\label{sec5}

We present the Proof of Theorem \ref{thlambda}.
\proof The claim $(b)$ follows immediately from $(a)$ and Remark \ref{remeq}. Indeed, in view
of Remark \ref{remeq}, a solution of \fbi\, with $\ga_{\sscp I}\geq 0$ exists if and only if
a solution of $\prl$ exists with $I=\lm^q$. It then follows from $(a)$ that no
solutions of \fbi\, exists with $\ga_{\sscp I}\geq 0$ whose
$I$ is larger than $\ov{\lm}^{\frac{p}{p-1}}$.\\
Therefore we are left with the proof of $(a)$. We argue by contradiction and assume that there exists
a sequence of solutions $(\al_n,\psi_n)$ of $\left.\prl\right|_{\lm=\lm_n}$ such that
$\lm_n\to +\ii$.
We will need the following,
\ble\label{minen}
Let $\psi$ be any a solution of
$$
\graf{-\Delta \psi=f \quad \mbox{in}\;\;\om\\ \psi=0 \quad \mbox{on}\;\;\pa \om}
$$
where $\ino f=1$ and $\ino |f|^N\leq C$. Then,
$$
\ino|\nabla \psi|^2\geq c>0,
$$
for some positive constant $c>0$ depending only by $C$, $N$ and $\om$.
\ele
\proof
Since $f\in L^{N}(\om)$ then by standard elliptic regularity theory
$\psi\in W^{2,N}_0(\om)$ and in particular we have
$$
\ino  |\nabla \psi|^2=\ino f \psi=\ino fG[f],
$$
for any such $\psi$ and $f$.
Therefore we are left to prove that
$$
\inf \left\{\ino fG[f]\,|\, \ino |f|^N\leq C,\,\ino f=1\right\}\geq  c>0.
$$

We argue by contradiction. If the claim were false we could find a sequence
$f_n$ such that $\ino |f_n|^N\leq C$, $\ino f_n=1$ and $\ino f_n G[f_n]\to 0^+$.
Since $f_n$ is uniformly bounded in $L^{N}(\om)$ we can pass to a subsequence such that
$f_n\rightharpoonup f_{\ii}$ weakly in $L^N(\om)$. Therefore, in particular we have
$\ino f_{\ii}=1$. Moreover, since $G[f_n]$ is uniformly bounded in $W^{2,N}_0(\om)$, then
by the Sobolev embedding we can pass to a further subsequence (which we will not relabel)
such that $G[f_n]\to \psi_{\ii}$ in
$C^{0}(\ov{\om})\cap W^{1,2}_0(\om)$. On the other
side, since for fixed $x\in\om$ in particular we have $G(x,\cdot)\in L^{\frac{N}{N-1}}(\om)$
(\cite{stam}), and
since $f_n\rightharpoonup f_{\ii}$ weakly in $L^N(\om)$,
then we see that $G[f_n]\to G[f_\ii]$ pointwise in $\om$. Therefore $G[f_\ii]\equiv\psi_\ii$ and in particular we conclude that
$$
0=\lim\limits_{n}\ino f_n G[f_n]=\ino f_\ii G[f_\ii],\quad f_\ii\in L^{N}(\om), \quad \ino f_{\ii}=1.
$$
This is clearly impossible since it is easy to check that then
$\psi_\ii \equiv \hspace{-.4cm} /\hspace{.2cm} 0$, $\psi_\ii\in W^{1,2}_0(\om)$
would satisfy $\ino |\nabla \psi_\ii|^2=\ino f_\ii G[f_\ii]=0$.

\finedim

\bigskip
By Lemma \ref{lemE1}, for any fixed $n$ it holds $\|\psi_n\|_{\ii}\leq C_n$.
Let $m_n=\sup\limits_{\om}(\al_n+\lm_n\psi_n)$. If there exists $C>0$ such that
$m_n\leq C$, then $\psi_n\leq \frac{C}{\lm_n}$ and consequently

$$
\ino |\nabla \psi_n|^2=\ino (\al_n+\lm_n\psi_n)^p\psi_n\leq \frac{C}{\lm_n}\to 0,\,\ainf.
$$

This contradicts Lemma \ref{minen} since $f_n=(\al_n+\lm_n\psi_n)^p$ would obviously satisfy the needed
assumptions in this case. Therefore we deduce that, passing to a subsequence if necessary,
$m_n\to +\ii$.
Let $x_n\in \om:m_n=\al_n+\lm_n\psi_n(x_n)$, we have only two possibilities:
either there exists a subsequence $x_n$ such that $x_n\to \ov{x}\in \om$, or dist$(x_n, \pa\om)\to 0$.
The second alternative can be ruled out by a well known moving plane argument based on \cite[Theorem 2.$1^{'}$]{gnn} and the Kelvin transform as in \cite[p. 52]{DLN}. Concerning the first
alternative without loss of generality we assume $x_n=0$ and define
$\dt_n=\frac{1}{\lm_n m_n^{p-1}}\to 0$, and
$$
v_n(y)=m_n^{-1}{(\al_n+\lm_n\psi_n(\dt_n y))},\quad y\in \om_n=(\dt_n)^{-1}\om
$$
which satisfies
$$
\graf{-\Delta v_n= v_n^p \quad \mbox{in}\;\;\om_n\\ \inf\limits_{\om_n}v_n\geq m_n^{-1}\al_n \\
\sup\limits_{\om_n}v_n\leq v_n(0)=1}
$$
Clearly for any $R\geq 1$, we have $R<\frac{\mbox{dist}(0, \pa\om_n)}{\dt_n}$, for $n$ large enough.
Therefore, passing to a subsequence if necessary, $v_n$ converges uniformly on compact subsets
of $\R^N$ to a solution of
$$
\graf{-\Delta v= v^p \quad \mbox{in}\;\;\R^N\\
v\geq 0 \quad \mbox{in}\;\;\R^N\\
\sup\limits_{\R^N}v\leq v(0)=1}
$$
Therefore in particular $v\in C^{2}(\R^N)$. This is clearly impossible since
it is well known (\cite{gs}) that the unique $C^{2}(\R^N)$, $N\geq 3$,
non negative solution of $-\Delta v=v^p$ in $\R^N$ is $v\equiv 0$. The same holds
for $N=2$ as it is well known that a superharmonic function which is also bounded from below in $\R^2$
must be constant, which is the desired contradiction.
\finedim

\bigskip
\bigskip

\section{\bf Generic properties}\label{sec6}
In this section we prove Theorems \ref{alt} and \ref{lmconvex}. As a first step,
by an adaptation of some arguments in \cite{H}, we will prove the following,
\bte\label{thm1.generic}
For any $\om_0\subset \R^N$ of class $C^{4}$ there exists a meager set
$\mathcal{F}\subset \mbox{\rm Diff}^{\,4}(\om_0)$, depending also on $N$ and $p$,
such that if $h\in \mbox{\rm Diff}^{\,4}(\om_0)\setminus\mathcal{F}$ then, for any positive
solution $(\all,\pl)$ of \mbox{\rm $\prl$} on $\om:=h(\om_0)$ with $\lm>0$, it holds: either\\
$(a)$   Ker$(L_{\lm})=\emptyset$, or\\
$(b)$  Ker$(L_{\lm})=\mbox{\rm span}\{\phi\}$ is one dimensional and
$<\phi>_{\ssl}\neq 0$.
\ete
Let us recall few definitions and set some notations first.

\bdf\label{C1} {\it
A domain $\om$ is of class $C^{k}(C^{k,r})$, $k\geq 1$, if for each $x_0\in \pa\om$ there exists a ball $B=B_r(x_0)$
and a one to one map $\Theta: B\mapsto U\subset \R^2$ such that $\Theta\in C^{k}(B)(C^{k,r}(B)),
\Theta^{-1}\in C^{k}(U)(C^{k,r}(U))$ and the following holds:
$$
\Theta(\om\cap B)\subset \R^2_+\quad \mbox{ and }\quad \Theta(\om\cap B)\subset \pa\R^2_+.
$$

It is well known {\rm(}see for example \cite{H}{\rm)} that this is equivalent to say that there exists
$r>0$ and $M>0$ such that, given any ball
$B_r(x_0)$, $x_0\in \R^2$ then, after suitable rotation and translations, it holds:
$$
\om\cap B=\{(x_1,x_2)\,:\,x_2<f(x_1)\}\cap B\quad \mbox{ and } \quad \pa\om\cap B=\{(x_1,x_2)\,:\,x_2=f(x_1)\}\cap B,
$$

for some $f\in C^{k}(\R)(C^{k,r}(\R))$ whose norm is not larger than $M.$}
\edf

\bdf {\it
Let $\om\subset \R^2$ be an open and bounded domain of class $C^m$, $m\geq 1$.
$C^m(\ov{\om}\,;\R^2)$, is the Banach space of continuous and $m$-times differentiable maps on $\om$,
whose derivatives of order $j=0,1,\cdots,m$
extend continuously on $\ov{\om}$. $\mbox{\rm Diff}^{\,m}(\om)\subset C^m(\ov{\om}\,;\R^2)$ is the open subset of
$C^m(\ov{\om}\,;\R^2)$
whose elements are $C^m$ imbeddings on $\ov{\om}$, that is, of maps $h:\ov{\om}\mapsto \R^2$ which are diffeomorphisms of class $C^m$
on their images ${h(\ov{\om}\,)}$.}
\edf

We recall that if $X,Z$ are Banach spaces and $T:X\to Z$ is linear and continuous, then T is Fredholm (semi-Fredholm)
if $R(T)$ (the range of $T$) is closed and both (at least one of) dim(Ker$(T)$) and codim$(R(T))$ are finite. If $T$ is Fredholm, then
the index of $T$ is
$$
\mbox{\rm ind}(T)=\mbox{\rm dim(Ker$(T)$)}-\mbox{\rm codim}(R(T)).
$$
Also, a semi-Fredholm operator with finite index is a Fredholm operator. We refer to \cite{Ka} for details and proofs.
Given a Banach space
$X$ and $x\in X$, we will denote by $T_x X$ the tangent space at $x$.

\bdf{\it
Let $X,Z$ be Banach spaces, $A\subset X$ an open set and $F:A\to Z$ a $C^1$ map. Suppose that for any $x \in  A$ the
Fr\'echet derivative $D_x F(x):T_xX\to T_\eta Z$ is a Fredholm operator. A point $x\in A$ is a \un{regular point}
if $D_x F(x)$ is surjective, is a \un{singular point} otherwise. The image of a singular point $\eta=F(x)\in Z$ is a \un{singular value}.
The complement of the set of singular values in $Z$ is the set of \un{regular values}.}
\edf

The following Theorem is a particular case of a more general transversality result proved in \cite{H}, see also \cite{Sm}.

\bte[{\rm \cite{H}}]\label{enry}
Let ${X},\mathcal{H},Z$ be separable Banach spaces, $\mathcal{A}\subseteq X\times \mathcal{H}$ an open set,\\
$\widetilde{\Phi}:\mathcal{A}\to Z$ a map of class $C^k$ and
$\eta\in Z$. Suppose that for each $(x,h)\in \widetilde{\Phi}^{-1}(\eta)$ it holds:
\begin{align}\nonumber
&(i) D_{x}\widetilde{\Phi}(x,h):T_x X\to T_\eta Z\;\mbox{is a Fredholm operator with index $< k$};\\
&(ii) D \widetilde{\Phi}(x,h)=(D_{x}\widetilde{\Phi}(x,h),D_{h}\widetilde{\Phi}(x,h)): T_x X\times T_h\mathcal{H} \to T_\eta Z\;\mbox{is surjective}.\nonumber
\end{align}
Let $A_h=\{x\,:\,(x,h)\in \mathcal{A}\}$ and
$$
\mathcal{H}_{\rm crit}=\{h\,:\,\eta \mbox{ is a singular value of }\widetilde{\Phi}(\,\cdot\,, h):A_h\to Z\}.
$$
Then $\mathcal{H}_{\rm crit}$ is meager in $\mathcal{H}$.
\ete

\bigskip
\bigskip

We are ready to present the proof of Theorem \ref{thm1.generic}.\\
{\it The Proof of Theorem \ref{thm1.generic}}\\
Let $\om_0$ as in the statement and let us define
$$
X_{\sscp \om_0}=\R\times\R \times C^{2,r}_0(\ov{\om_0}\,).
$$
As in section \ref{sec2}, see Remark \ref{Analytic}, we denote by $A_{\sscp \om_0}$ an open subset of the Banach space of triples
$(\lm,\al,\psi)\in X_{\sscp \om_0}$ such that, on $A_{\sscp \om_0}$, the density $\rla=\rla(\psi)=(\al+\lm\psi)^p$ is well defined and
$$
\all+\lm\pl\geq \frac{\all}{2} \quad \mbox{ in }\ov{\om_0}
$$
in a sufficiently small open neighborhood in $A_{\sscp \om_0}$ of
any triple of the form $(\lm,\all,\pl)$ whenever $(\all,\pl)$ is a positive solution of {\rm $\prl$}.
In what follows we write
$$
\rla=\rla(\psi)=(\al+\lm \psi)^{p},\quad (\lm,\al,\psi)\in A_{\sscp \om_0}.
$$

We define the maps,
$$
F_{\om_0}:A_{\sscp \om_0} \to C^r(\ov{\om_0}\,  ),
\quad  F_{\om_0}(\lm,\al,\psi)=\Delta \psi+ \rla(\psi),
$$
$$
{M}_{\om_0}:A_{\sscp \om_0}\to \R,
\quad {M}_{\om_0}(\lm,\al,\psi)=\left(\,\int\limits_{\om_0} \rla(\psi)\right)-1.
$$

Next, for fixed $h\in \mbox{\rm Diff}^{\, 4}(\om_0)$ and
$\psi\in C^{2,r}_0(\ov{h(\om_0)}\,)$,
we define the pull back,
$$
h^*(\psi)(x)=\psi(h(x)),\;x\in \ov{\om_0}.
$$
Clearly $h^*$ is an isomorphism of $C^{2,r}_0(\ov{h(\om_0)}\,)$ onto $C^{2,r}_0(\ov{\om_0}\,)$
with inverse $h^{*-1}=(h^{-1})^*$. For any such $h$, it is well defined the map
$$
F_{h(\om_0)}:A_{h(\om_0)}\to C^{r}(\ov{h(\om_0)}\,)
$$
and then we can set,
$$
h^*F_{h(\om_0)}h^{*-1}:A_{\om_0}\times \mbox{\rm Diff}^{\, 4}(\om_0)\to C^{r}(\ov{\om_0}\,),
$$
and
$$
M^*_{\om_0}:A_{\om_0}\times \mbox{\rm Diff}^{\, 4}(\om_0)\to \R,
$$
where $M^*$ is defined as follows,
$$
M^*_{\om_0}(\lm,\al,\psi,h)=M_{h(\om_0)}(\lm,\al,(h^*)^{-1}(\psi)).
$$

Putting $\mathcal{H}=\mbox{\rm Diff}^{\, 4}(\om_0)$,
$\eta=0\in Z=C^{r}(\ov{\om_0}\,)$, we will apply Theorem \ref{enry} to the map $\widetilde{\Phi}=\widetilde{\Phi}(\lm,\al,\psi,h)$ defined as follows

$$
\widetilde{\Phi}:\mathcal{A}\to \R\times C^{r}(\ov{\om_0}\,),\qquad
\mathcal{A}=A_{\sscp \om_0}\times \mathcal{H},
$$

$$
\widetilde{\Phi}(\lm,\al,\psi,h)=\left( \begin{array}{c}\widetilde{\Phi}_1(\lm,\al,\psi,h) \\ $\,$\\
\widetilde{\Phi}_2(\lm,\al,\psi,h)\end{array}\right)=\left( \begin{array}{c}h^*F_{h(\om_0)}h^{*-1}(\lm,\al,\psi) \\ $\,$\\
{M}^*_{\om_0}(\lm,\al,\psi,h)\end{array}\right).
$$

\bigskip
\bigskip
{\bf STEP 1:} Our aim is to show that the assumptions $(i)$ and $(ii)$ of Theorem \ref{enry} hold.\\
As in \cite{H}, it is very useful for the discussion to denote by
$(\dot{\lm},\dot{\al},\dot{\psi},\dot{h})\in \R\times \R\times C^{2,r}_0(\ov{\om_0}\,)\times C^{4}(\ov{\om_0}\,;\,\R^2)$ the
elements of the tangent space at points
$({\lm},{\al},{\psi},{h})\in A_{\sscp \om_0} \times \mathcal{H}$.\\
First of all observe that for fixed $h\in \mbox{\rm Diff}^{\, 4}(\om_0)$,
the linearized  operator,
$$
D_{\lm,\al,\psi}\widetilde{\Phi}_1(\lm,\al,\psi,h):\R\times \R\times C^{2,r}_0(\ov{\om_0}\,) \to C^{r}(\ov{\om_0}\,),
$$
acts as follows on a triple $(\dot{\lm},\dot{\al},\dot{\psi})\in \R\times \R\times C^{2,r}_0(\ov{\om_0}\,)$,
$$
D_{\lm,\al,\psi}\widetilde{\Phi}_1(\lm,\al,\psi,h)[\dot{\lm},\dot{\al},\dot{\psi}]=
h^*\left(\Delta \dot{\psi}^*+(\rla(\psi^*))^{\frac{1}{q}}(\lm p \dot{\psi}^*+p{\psi}^*\dot{\lm}
+p\dot{\al})\right),
$$
where
$$
{\psi}^*=(h^{*})^{-1}{\psi},\quad \dot{\psi}^*=(h^{*})^{-1}\dot{\psi}.
$$
Since any diffeomorphism of class $C^{4}$ maps the Laplace operator to a uniformly elliptic operator
with $C^2$ coefficients, by standard elliptic estimates it is not difficult  to see that
$D_{\lm,\al,\psi}\widetilde{\Phi}_1(\lm,\al,\psi,h)$ is a Fredholm operator
of index 2. Moreover we have,
$$
D_{\lm,\al,\psi}\widetilde{\Phi}_2(\lm,\al,\psi,h):\R\times \R\times C^{2,r}_0(\ov{\om_0}\,) \to \R,
$$
which acts as follows,
$$
D_{\lm,\al,\psi}\widetilde{\Phi}_2(\lm,\al,\psi,h)[\dot{\lm},\dot{\al},\dot{\psi}]=
\int\limits_{h(\om_0)}(\rla(\psi^*))^{\frac{1}{q}}(\lm p\dot{\psi}^*+p{\psi}^*\dot{\lm}+p\dot{\al})
$$

and is semi-Fredholm with index $+\ii$ and range of vanishing codimension. Therefore the codimension of the range
of $D_{\lm,\al,\psi}\widetilde{\Phi}(\lm,\al,\psi,h)$ is the same as that of $D\widetilde{\Phi}_1(\psi,\lm,\al,h)$. In particular,
since Ker$(D_{\lm,\al,\psi}\widetilde{\Phi}(\lm,\al,\psi,h))=$
Ker$(D_{\lm,\al,\psi}\widetilde{\Phi}_1(\lm,\al,\psi,h))\cap$Ker$(D_{\lm,\al,\psi}\widetilde{\Phi}_2(\lm,\al,\psi,h))$,
then we conclude that $D_{\lm,\al,\psi}\Phi(\lm,\al,\psi,h)$ is a Fredholm
operator of index at most $2$.\\
This fact proves $(i)$ whenever we can show that $\widetilde{\Phi}\in C^{k}(\mathcal{A})$ for some
$k\geq 3$.  The regularity of $\widetilde{\Phi}_1$ with respect to $h$ is the same as that of $F_{h(\om_0)}$ with respect to $\psi$,
see chapter 2 in \cite{H}. Therefore, in view of Remark \ref{Analytic}, we have $\widetilde{\Phi}_1\in C^{\ii}(\mathcal{A})$.
The regularity of $\widetilde{\Phi}_2$ is more subtle since derivatives of order $m\geq 2$ with respect to $h$
also involve the derivatives of the normal field on $\pa\om$. However, since $\om$ is of class $C^4$ and
still by Remark \ref{Analytic}, we can apply Theorem 1.11 in \cite{H} to conclude that $\widetilde{\Phi}_2$ is of class $C^3$
with respect to $h$.
In particular $\widetilde{\Phi}_2\in C^3(\mathcal{A})$ and then also $\widetilde{\Phi}\in C^3(\mathcal{A})$, as claimed.\\ \\
Next we prove $(ii)$, that is, we show that $\eta=0$ is a regular value for the map $(\lm,\al,\psi,h)\to \widetilde{\Phi}(\lm,\al,\psi,h)$.
We argue by contradiction and
suppose that there exists a singular point $(\ov{\lm},\ov{\al},\ov{\psi},\ov{h}\,)$ of $\widetilde{\Phi}$ such that
$\widetilde{\Phi}_i(\ov{\lm},\ov{\al},\ov{\psi},\ov{h}\,)=0$,
$i=1,2$.\\
First of all, let us define $\om=\ov{h}(\om_0)$, $\ov{\varphi}=(\ov{h}^*)^{-1}\ov{\psi}\in C^{2,r}_0(\ov{\om}\,)$ and
$\widehat{\Phi}(\lm,\al,\varphi,h)$
on $A_{\sscp \om}\times \mbox{\rm Diff}^{\, 4}(\om)$ as follows,
$$
\widehat{\Phi}(\lm,\al,\varphi,h)=\left(\begin{array}{c}
                    \widehat{\Phi}_1(\lm,\al,\varphi,h) \\
                    \widehat{\Phi}_2(\lm,\al,\varphi,h)
                  \end{array}\right)=
\left( \begin{array}{c}h^*F_{\om}h^{*-1}(\lm,\al,\varphi) \\ $\,$\\
{M}^*_{\om}(\lm,\al,\varphi,h)\end{array}\right),
$$
where,
$$
F_{\om}:A_{\sscp \om}\to C^r(\ov{\om}),
\quad  F_{\om}(\lm,\al,\varphi)=\Delta \varphi+ \rla(\varphi),
$$
$$
{M}^*_{\om}:A_{\sscp \om}\to \R,
\quad{M}^*_{\om}(\lm,\al,\varphi,h)=M_{h(\om)}(\lm,\al,(h^*)^{-1}(\varphi)).
$$
Let $i_\om\in  \mbox{\rm Diff}^{\, 4}(\om)$ be the identity map.
By construction, in these new coordinates the map
$\widehat{\Phi}(\lm,\al,\varphi,h)$ has a singular point $(\ov{\lm},\ov{\al},\ov{\varphi},i_\om)$ such that
$\widehat{\Phi}_i(\ov{\lm},\ov{\al},\ov{\varphi},i_\om)=0$, $i=1,2$, that is,
by assumption the derivative $D_{\lm,\al,\varphi,h}\widehat{\Phi}(\ov{\lm},\ov{\al},\ov{\varphi},i_\om)$ is not surjective.
Putting
$$
\ov{\rh}=(\ov{\al}+\ov{\lm}\ov{\varphi})^p,
$$
a subtle evaluation shows that
$D_{\lm,\al,\varphi,h}\widehat{\Phi}_1(\ov{\lm},\ov{\al},\ov{\varphi}, i_\om)$ acts on
$$
(\dot{\lm},\dot{\al},\dot{\varphi},\dot{h})\in \R\times\R\times C^{2,r}_0(\ov{\om})\times C^{4}(\ov{\om}\,;\R^2)
$$
as follows (see Theorem 2.2 in \cite{H}),
\begin{align}
&\qquad D_{\lm,\al,\varphi,h}\widehat{\Phi}_1(\ov{\lm},\ov{\al},\ov{\varphi}, i_\om)[\dot{\lm},\dot{\al},\dot{\varphi},\dot{h}]\nonumber\\
&=\Delta \dot{\varphi}+p\left(\ov{\rh}\right)^{\frac1q}\left(\ov{\lm}\dot{\varphi}+\ov{\varphi}\dot{\lm}+ \dot{\al}\right)+\dot{h}\cdot \nabla(\Delta \ov{\varphi}+\ov{\rh})-(\Delta +p\ov{\lm}(\ov{\rh}\,)^{\frac1q})\dot{h}\cdot \nabla \ov{\varphi}\nonumber\\
&=\left(\Delta +p\ov{\lm}(\ov{\rh})^{\frac1q}\right)\dot{\varphi}-
\left(\Delta +p\ov{\lm}(\ov{\rh}\,)^{\frac1q}\right)\dot{h}\cdot \nabla \ov{\varphi}+p(\ov{\rh})^{\frac1q}( \ov{\varphi}\dot{\lm}+\dot{\al})\label{DF1},
\end{align}
where we used the fact that
$\Delta \ov{\varphi}+\ov{\rh}=\widehat{\Phi}_1(\ov{\lm},\ov{\al},\ov{\varphi}, i_\om)=0$.
In particular, see Theorem 1.11 in \cite{H}, we have,
$$
D_{\lm,\al,\varphi,h}\widehat{\Phi}_2(\ov{\lm},\ov{\al},\ov{\varphi}, i_\om)[\dot{\lm},\dot{\al},\dot{\varphi},\dot{h}]=
\ino p(\ov{\rh}\,)^{\frac1q}\left(\ov{\lm}\dot{\varphi}+\ov{\varphi}\dot{\lm}+\dot{\al}\right)+
\int\limits_{\pa\om}\ov{\rh}\dot{h} \cdot \nu
$$
At this point observe that, by the Fredholm property of the operator $\Delta +\ov{\rh}$ on $C^{2,r}_0(\ov{\om}\,)$,
we have that the subspace
$\left\{D_{\lm,\al,\varphi,h}\widehat{\Phi}_1(\ov{\lm},\ov{\al},\ov{\varphi}, i_\om)[(0,0,\dot{\varphi},0)],
\,\dot{\varphi}\in C^{2,r}_0(\ov{\om})\right\}$, is closed and has finite codimension. Next,
since $\ov{\varphi}\in C^{2,r}_0(\ov{\om}\,)$ and $\pa\om$ is of class $C^{4}$, then
by standard elliptic regularity theory we find that $\ov{\varphi}\in C^{3,r}_0(\ov{\om})$ and then
$\dot{h}\cdot \nabla \ov{\varphi}\in C^{2,r}(\ov{\om})$. As a consequence we can prove
that  the subspace
$\left\{ D_{\lm,\al,\varphi,h}\widehat{\Phi}_1(\ov{\lm},\ov{\al},\ov{\varphi}, i_\om)[(0,0,0,\dot{h})],
\dot{h}\in C^{4}(\ov{\om}\,;\R^2)\right\}$ is closed with finite codimension as well.
Indeed, let us define
$K:C^{2,r}(\ov{\om}\,)\mapsto C^{2,r}(\ov{\om}\,)$ as the
linear operator which, to any $\phi\in C^{2,r}(\ov{\om}\,)$, associates the unique solution
$\phi_b=K[\phi]$ of $\Delta \phi_b=0$,
$\phi_b=\phi$ on $\pa\om$.
Clearly this is always well posed since $\om$ is of class $C^4$ and $\phi\in C^{2,r}(\ov{\om}\,)$. Then,
$$
\Delta \phi + p\ov{\lm}(\ov{\rh}\,)^{\frac1q}\phi=g\in C^{r}(\ov{\om}\,),
$$
if and only if
$$
\phi \in C^{2,r}(\ov{\om}\,)\mbox{ and }\phi+T[\phi]=G[g]\in C^{2,r}(\ov{\om}\,),
$$
where $G[g]=\ino G(x,y)g(y)$ and $T:C^{2,r}(\ov{\om}\,)\mapsto C^{2,r}(\ov{\om}\,)$,
$T(\phi)=G[p\ov{\lm}(\ov{\rh}\,)^{\frac1q}\phi]-K[\phi]$.
Since $\om$ is of class $C^{4}$, then by standard elliptic estimates (\cite{GT}), $T$
maps $C^{2,r}(\ov{\om}\,)$ into $C^{3,r}(\ov{\om}\,)$. Therefore $T$ is compact and
then we conclude by the Fredholm alternative that the range of
$(\Delta + p\ov{\lm}(\ov{\rh}\,)^{\frac1q})(\dot{h}\cdot \nabla \ov{\varphi})$,
$\dot{h}\cdot \nabla \ov{\varphi}\in C^{2,r}(\ov{\om}\,)$,
is closed in $C^{r}(\ov{\om}\,)$ and has finite codimension.\\

At this point we readily deduce from these two facts that there exists
$(t_{\perp},\phi_{\perp})\in \R\times C^r(\ov{\om}\,)$ which is orthogonal to the image of
$D_{\lm,\al,\varphi,h}\widehat{\Phi}(\ov{\lm},\ov{\al},\ov{\varphi}, i_\om)$, that is,
\beq\label{1perp}
\ino \phi_{\perp}
\left(\left(\Delta +p\ov{\lm}(\ov{\rh}\,)^{\frac1q}\right)\dot{\varphi}-
\left(\Delta +p\ov{\lm}(\ov{\rh}\,)^{\frac1q}\right)\dot{h}\cdot \nabla \ov{\varphi}+p(\ov{\rh}\,)^{\frac1q}( \ov{\varphi}\dot{\lm}+\dot{\al})\right)=0,
\forall\,(\dot{\lm},\dot{\al},\dot{\varphi},\dot{h}),
\eeq
and
$$
t_{\perp}\left(\ino p(\ov{\rh}\,)^{\frac1q}\left(\ov{\lm}\dot{\varphi}+\ov{\varphi}\dot{\lm}+\dot{\al}\right)+
\int\limits_{\pa\om}\ov{\rh}\dot{h} \cdot \nu\right)=0,\forall\,(\dot{\lm},\dot{\al},\dot{\varphi},\dot{h}).
$$
Putting $(\dot{\lm},\dot{\al},\dot{\varphi},\dot{h})=(0,\dot{\al},0,0)$ in the second equation we see that
$t_{\perp}=0$. Next, putting $(\dot{\lm},\dot{\al},\dot{\varphi},\dot{h})=(0,\dot{\al},0,0)$ in \rife{1perp}
we find $\ino (\ov{\rh}\,)^{\frac1q}\phi_{\perp}=0$, while
if we choose $(\dot{\lm},\dot{\al},\dot{\varphi},\dot{h})=(\dot{\lm},0,0,0)$ we have $\ino (\ov{\rh}\,)^{\frac1q}\ov{\varphi}\phi_{\perp}=0$.
At this point, putting $\dot{h}=0$, we find that,
$$
\ino \phi_{\perp}\left(\Delta +p\ov{\lm}(\ov{\rh}\,)^{\frac1q}\right)\dot{\varphi}=0,\,\forall \, \dot{\varphi}\in C^{2,r}_0(\ov{\om}\,),
$$
which shows that $\phi_{\perp}$ is a $C^{r}_0(\ov{\om}\,)$ solution of $\Delta  \phi_{\perp}+p\ov{\lm}(\ov{\rh}\,)^{\frac1q}\phi_{\perp}=0$ in the $(C^2_0(\ov{\om}))^*$ sense.
Therefore, by standard elliptic estimates (where we recall that $\pa\om$ is of class $C^{4}$),
$\phi_{\perp}$ is a $C_0^2(\ov{\om}\,)$ solution of $\Delta  \phi_{\perp}+p\ov{\lm}(\ov{\rh}\,)^{\frac1q}\phi_{\perp}=0$.
As a consequence we observe that \rife{1perp} is reduced to
$$
\ino \phi_{\perp}\left(\Delta +p\ov{\lm}(\ov{\rh}\,)^{\frac1q}\right)\dot{h}\cdot \nabla\ov{\varphi}=0,\quad \forall\, \dot{h}\in C^{4}(\ov{\om}\,;\R^2),
$$
which allows us to deduce that,
\begin{align*}
0 &=\ino \phi_{\perp}\left(\Delta +p\ov{\lm}(\ov{\rh}\,)^{\frac1q}\right)\dot{h}\cdot \nabla\ov{\varphi}\nonumber\\
&=\ino \phi_{\perp}\left(\Delta +p\ov{\lm}(\ov{\rh}\,)^{\frac1q}\right)\dot{h}\cdot \nabla \ov{\varphi}-
\ino \left(\Delta \phi_{\perp}+p\ov{\lm}(\ov{\rh}\,)^{\frac1q}\phi_{\perp}\right)\dot{h}\cdot \nabla \ov{\varphi}\nonumber\\
&=\ino \phi_{\perp}\Delta(\dot{h}\cdot \nabla \ov{\varphi})-\ino (\Delta\phi_{\perp})\dot{h}\cdot \nabla \ov{\varphi}=\int\limits_{\pa\om} \left(\phi_{\perp}\pa_{\nu}(\dot{h}\cdot \nabla \ov{\varphi})-\dot{h}\cdot \nabla \ov{\varphi}(\pa_\nu\phi_{\perp})\right)\nonumber\\
&=-\int\limits_{\pa\om} (\pa_\nu\phi_{\perp})\dot{h}\cdot \nabla \ov{\varphi}\nonumber=-\int\limits_{\pa\om} (\pa_\nu\phi_{\perp})(\pa_\nu \ov{\varphi})\dot{h}\cdot \nu,\quad \forall\,\dot{h}\in C^4(\ov{\om},\R^2).
\end{align*}

Therefore, since $\dot{h}$ is arbitrary, we conclude that,
$$
(\pa_\nu\phi_{\perp})(\pa_\nu \ov{\varphi})\equiv0 \quad \mbox{ on }\pa\om.
$$
At this point we observe that since $\ov{\rh}>0$  on $\ov{\om}$ and $\ov{\varphi}=0$ on $\pa\om$, then, by the strong maximum principle,
we have $\ov{\varphi}>0$ in $\om$. Since $\pa\om$ is of class $C^{4}$ we can apply the Hopf boundary Lemma and conclude that
$\pa_\nu \ov{\varphi}<0$ on $\pa\om$. Therefore we conclude that necessarily $\pa_\nu\phi_{\perp}\equiv 0$ on $\pa\om$, which is in contradiction
with the Hopf boundary Lemma. This contradiction shows that $(ii)$ holds and then we can apply Theorem \ref{enry} and conclude that
there exists a meager set $\mathcal{F}\subset \mbox{Diff}^4(\om_0)$ such that if $h(\om_0)\notin \mathcal{F}$ then $\eta=0$ is a regular
value of $\widetilde{\Phi}(\lm,\al,\psi,h)$.\\

{\bf STEP 2:}
We have from STEP 1 that there exists a meager set $\mathcal{F}\subset \mbox{Diff}^4(\om_0)$ such that if $\om:=h(\om_0)\notin \mathcal{F}$,
then $\eta=0$ is a regular value of the map $\widetilde{\Phi}(\lm,\al,\psi,h)$.
As a consequence, for any $(\ov{\lm},\ov{\al},\ov{\psi})$ which solves

$$
\Phi(\lm,\al,\psi)=\left( \begin{array}{c}\Phi_1(\lm,\al,\psi)\\ $\,$\\
\Phi_2(\lm,\al,\psi)\end{array}\right)=\left( \begin{array}{c}F_{\om}(\lm,\al,\psi) \\ $\,$\\
-1+\ino \rla(\psi)\end{array}\right)=\left( \begin{array}{c}0 \\ $\,$\\
0\end{array}\right)
$$
and setting $\ov{\rh}=(\ov{\al}+\ov{\lm}\,\ov{\psi})^p$, then the differential
 $$
\left(\begin{array}{c}
                        D_{\lm,\al,\psi}\Phi_1(\ov{\lm},\ov{\al},\ov{\psi})[\dot{\lm},\dot{\al},\dot{\psi}]\\
                         D_{\lm,\al,\psi}\Phi_2(\ov{\lm},\ov{\al},\ov{\psi})[\dot{\lm},\dot{\al},\dot{\psi}]
                      \end{array}\right)=\left(\begin{array}{c}
                        \Delta\dot{\psi}  +p\ov{\lm}(\ov{\rh}\,)^{\frac1q}\dot{\psi} +
                        p(\ov{\rh}\,)^{\frac1q}(\ov{\psi}\dot{\lm}+\dot{\al})\\
                         \ino p(\ov{\rh}\,)^{\frac1q}(\ov{\lm}\dot{\psi} +\ov{\psi}\dot{\lm}+\dot{\al})
                      \end{array}\right)
                      $$
is surjective. In particular, there exists a subspace $V$ of the kernel of
$D_{\lm,\al,\psi}\Phi_2(\ov{\lm},\ov{\al},\ov{\psi})$,
$$
V\subseteq \mbox{Ker}\left(D_{\lm,\al,\psi}\Phi_2(\ov{\lm},\ov{\al},\ov{\psi})\right),
$$
such that $D_{\lm,\al,\psi}\Phi_1(\ov{\lm},\ov{\al},\ov{\psi})$ is surjective along $V$.
Let $\widetilde{V}$ be the {space} defined as follows
$$
\widetilde{V}=\{(\dot{\lm},\dot{\psi}):(\dot{\lm},\dot{\al},\dot{\psi})\in V\}
$$
Since $\dot{\al}+\ov{\lm}<\dot{\psi}>_{\sscp \ov{\lm}}+\dot{\lm}<\ov{\psi}>_{\sscp \ov{\lm}}=0$ in $V$,
then $\widetilde{V}$ cannot be {trivial}. Therefore we can write
$$
\dot{\al}=-\ov{\lm}<\dot{\psi}>_{\sscp \ov{\lm}}-\dot{\lm}<\ov{\psi}>_{\sscp \ov{\lm}},\quad \forall
(\dot{\lm},\dot{\psi})\in \widetilde{V},
$$
and substitute in $D_{\lm,\al,\psi}\Phi_1(\ov{\lm},\ov{\al},\ov{\psi})[\dot{\lm},\dot{\al},\dot{\psi}]$
to conclude that the map
$$
\widetilde{L}[\dot{\lm},\dot{\psi}]:=
\Delta\dot{\psi}  +p\ov{\lm}(\ov{\rh}\,)^{\frac1q}[\dot{\psi}]_{\sscp \ov{\lm}} +
p(\ov{\rh}\,)^{\frac1q}[\,\ov{\psi}\,]_{\sscp \ov{\lm}}\dot{\lm}
$$
with domain $(\dot{\lm},\dot{\psi})\in \widetilde{V}$ is surjective.\\
On the other side, since $(\ov{\al},\ov{\psi})$ is a solution of $\prl$,
then the operator,
$$
\ov{L}[\dot{\psi}]=\Delta\dot{\psi}  +p\ov{\lm}(\ov{\rh}\,)^{\frac1q}[\dot{\psi}]_{\sscp \ov{\lm}}
$$
is just \rife{eLl} evaluated at $(\lm,\all,\pl)=(\ov{\lm},\ov{\al},\ov{\psi})$
and we will use the fact that the Fredholm alternative holds for $\ov{L}$.\\
Let us define $R=R(\ov{L})\subseteq C^{r}(\ov{\om}\,)$  to be the range of $\ov{L}$.
Now if $\widetilde{V}=\{0\}\times W$ where $W$ is some vector subspace of
$C^{2,r}_0(\ov{\om}\,)$, then we have that $\ov{L}$ is surjective on $W$ and then
a fortiori $\ov{L}$ is surjective on $C^{2,r}_0(\ov{\om}\,)$.
Therefore, by the Fredholm alternative, we have Ker$(\ov{L})=\emptyset$, which is $(a)$ in the statement of
Theorem \ref{thm1.generic}. As a consequence we can assume that
$\widetilde{V}=\R\times W$, where $W$ is some vector subspace of
$C^{2,r}_0(\ov{\om}\,)$. As above, we can assume without loss of generality that
$\widetilde{L}$ is surjective on $\R\times C^{2,r}_0(\ov{\om}\,)$, that is, in other words that
$\widetilde{V}=\R\times C^{2,r}_0(\ov{\om}\,)$. At this point, let
$\ov{d}=\mbox{codim}(R)$ be the codimension of $R$.
Since $\widetilde{L}$ is surjective, then it is not difficult to see that $\ov{d}\leq 1$.
If $\ov{d}=0$, then
 Ker$(\ov{L})=\emptyset$ which is $(a)$ in the statement of Theorem \ref{thm1.generic}. Otherwise $\ov{d}=1$
 and we will conclude the proof by showing that $(b)$ holds in this case. In this situation the kernel must be one dimensional,
 Ker$(\ov{L})=\mbox{\rm span}\{\ov{\phi}\}$, for some $\ov{\phi}\in C^{2,r}_0(\ov{\om}\,)$ which
 satisfies $\Delta \ov{\phi}+p\ov{\lm}(\ov{\rh}\,)^{\frac1q}[\ov{\phi}]_{\sscp \ov{\lm}}=0$ and
 \rife{2907.6n} takes the form

 \beq\label{3}
  (\ov{\al}+p\ov{\lm}<\ov{\psi}>_{\sscp{\ov{\lm}}})<\ov{\phi}>_{\sscp{\ov{\lm}}}=\lm(p-1)<\ov{\psi}\,\ov{\phi}>_{\sscp{\ov{\lm}}}.
 \eeq
It is worth to recall that $\ov{\al}>0$, $\ov{\lm}>0$, $p>1$ by assumption and that,
by the maximum principle, $\ov{\psi}>0$ in $\om$.
Therefore it is enough to show that if $\ov{d}=1$ then $<\ov{\psi}\,\ov{\phi}>_{\sscp{\ov{\lm}}}\neq 0$.\\
We argue by contradiction and suppose that $<\ov{\psi}\,\ov{\phi}>_{\sscp{\ov{\lm}}}=0$.
Since  $\widetilde{L}$ is surjective, then $(\ov{\rh}\,)^{\frac1q}\ov{\phi}$
must be an element of its range and then there exists $(\mu,\phi)\in \R\times C^{2,r}_0(\ov{\om}\,)$ which
satisfies
$$
\Delta\phi  +p\ov{\lm}(\ov{\rh}\,)^{\frac1q}[\phi]_{\sscp \ov{\lm}} +
p(\ov{\rh}\,)^{\frac1q}[\,\ov{\psi}\,]_{\sscp \ov{\lm}}\mu=(\ov{\rh}\,)^{\frac1q}\ov{\phi}.
$$
Recalling that $<\ov{\phi}[\phi]_{\sscp \ov{\lm}}>_{\sscp \ov{\lm}}\,=<[\,\ov{\phi}\,]_{\sscp \ov{\lm}}\phi>_{\sscp \ov{\lm}}$,
multiplying this equation by $\ov{\phi}$ and integrating by parts we conclude that
$$
p\mu<\ov{\phi}[\,\ov{\psi}\,]_{\sscp \ov{\lm}}>_{\sscp \ov{\lm}}\,=<\ov{\phi}^2>_{\sscp \ov{\lm}}>0,
$$
which also shows that necessarily $\mu\neq 0$ and $<\ov{\phi}[\,\ov{\psi}\,]_{\sscp \ov{\lm}}>_{\sscp \ov{\lm}}\neq 0$.
However this is impossible since, in view of $<\ov{\psi}\,\ov{\phi}>_{\sscp{\ov{\lm}}}=0$ and \rife{3},
we would also conclude that
$$
0<p\mu <\ov{\phi}\,[\,\ov{\psi}\,]_{\sscp \ov{\lm}}>_{\sscp \ov{\lm}}=
p\mu<\ov{\phi}\,\ov{\psi}>_{\sscp{\ov{\lm}}}-p\mu<\ov{\phi}>_{\sscp{\ov{\lm}}}<\ov{\psi}>_{\sscp{\ov{\lm}}}=0,
$$
which is the desired contradiction.

\finedim

\bigskip
\bigskip

We are ready to present the proof of Theorem \ref{alt}.\\
{\it The Proof of Theorem \ref{alt}.}\\
By Theorem B we have that $\sg_1(\all,\pl)>0$ and $\all>0$ for any
$\lm\in[0,\frac1p\Lambda(\om,2p)]$. Therefore
by Lemma \ref{lem1.1} we can continue the curve $(\all,\pl)$ in a right neighborhood of $\frac1p\Lambda(\om,2p)$
to a larger simple real analytic curve
of solutions with no bifurcation points,
$\mathcal{G}(\om)\subset \mathcal{G}_{\mu}$, $\mu>\frac1p\Lambda(\om,2p)$, such that,  by continuity,
$\all>0$ and $\sg_1(\all,\pl)>0$ for any $\lm<\mu$.
Therefore it is well defined,
$$
{\lm}_1(\om):=\sup\left\{\mu>\frac1p\Lambda(\om,2p)\,:\,\all>0\mbox{ and }\sg_1(\all,\pl)>0,\;
\forall\,(\all,\pl)\in \mathcal{G}_{\mu},\,\forall \lm<\mu\right\}.
$$
By Theorem \ref{thlambda} we have $\lm_1(\om)<+\ii$ and there are only two possibilities:\\
either $\inf\limits_{\lm\in[0,\lm_1(\om))}\all=0$ or $\inf\limits_{\lm\in[0,\lm_1(\om))}\all>0$.\\
By Proposition \ref{pr-enrg} we have that
$\frac{d\all}{d\lm}<0,\frac{d\el}{d\lm}>0$ continue to hold whenever $\all>0$
and $\sg_1(\all,\pl)>0$ are satisfied,
whence for any $\lm<\lm_1(\om)$ as well. Now if $\inf\limits_{\lm\in[0,\lm_1(\om))}\all=0$,
then by definition we have that $\all>0$ and $\sg_1(\all,\pl)>0$ for any $\lm<\lm_1(\om)$, and then
by Lemma \ref{lemE1} it is not difficult to see that $(i)$ holds with $\lm_\ii(\om,p)=\lm_1(\om)$.\\

We conclude the proof by showing that if $\al_{\sscp \lm_1(\om)}>0$ and
$\om=|h(\om_0)|^{-\frac1N}h(\om_0)$, where
$h\in \mbox{\rm Diff}^{\,4}(\om_0)\setminus\mathcal{F}$ as in
the statement of Theorem \ref{thm1.generic}, then, possibly taking a larger but still meager set
$\mathcal{F}$, $(ii)$ holds.\\

Let $\widetilde{\om}=h(\om_0)$, then
either $(a)$ or $(b)$ of Theorem \ref{thm1.generic} hold for $\prl$ on $\widetilde{\om}=h(\om_0)$,
which however in general does not satisfy $|\widetilde{\om}|=1$. On the other side,
of course by the dilation invariance of $\prl$, $(ii)$ holds for $\om$ if and only if $(ii)$
holds on $\widetilde{\om}$. Therefore we can assume without loss of generality that
either $(a)$ or $(b)$ of Theorem \ref{thm1.generic} hold for $\prl$ on
${\om}$ which satisfies $|\om|=1$.
We need first the following,
\ble\label{altlem}
$\mathcal{G}_{\lm_1(\om)}(\om)$
can be continued to a continuous simple parametrization without bifurcation points,
$$
\mathbb{G}_\ii=\left\{[0,s_{\ii})\ni s\mapsto (\lm(s),\al(s),\psi(s))\right\}
$$
which has locally the structure of a $1$-dimensional real analytic manifold,
such that for any $s$, $(\al(s),\psi(s))$ is a positive solution of
{\rm $\left.\prl\right|_{\lm=\lm(s)}$} and, as $s\to s_{\ii}$, we have $\al(s)\to 0^+$,
and  for any sequence $s_n\to s_\ii$, there exists a subsequence $\{t_n\}\subseteq \{s_n\}$,
such that $\lm(t_n)\to \lm_\ii$ and
$\psi(t_n)\to \psi_{\ii}$ in $C_0^2(\ov{\om}\,)$, for some
$\lm_\ii\in (\frac1p\Lambda(\om,2p),\ov{\lm}(\om,p)]$
and $\psi_{\ii}$ which solves {\rm $\prl$} with $\lm=\lm_\ii$ and $\al=0$.
\ele
\proof By Theorem \ref{thlambda} and Lemma \ref{lemE1} we see that $\lm$ and $\pl$ are compact
in the $\R$ and $C^{2}_0(\ov{\om})$ topology respectively. In view also of Theorem B,
any limit point $\lm$ of a sequence of solutions whose $\al_n\to 0^+$ must satisfy
$\lm\in(\frac1p\Lambda(\om,2p),\ov{\lm}(\om,p)]$.
Therefore we just need to prove that $\mathcal{G}_{\lm_1(\om)}(\om)$
can be continued to $\mathbb{G}_\ii$ and that $\al(s)\to 0^+$ as $s\to s_{\ii}$.\\
By Lemmas \ref{lem1.1} and \ref{lemE1}, and since $\al_{\sscp \lm_1(\om)}>0$ by assumption,
it is not difficult to see that we must have $\sg_1(\al_1,\psi_1)=0$,
where $\al_1=\al_{\sscp \lm_1(\om)}>0$ and $\psi_1=\psi_{\sscp \lm_1(\om)}$. Thus, in view
of Theorem \ref{thm1.generic}, we have
$Ker(L_{\sscp \lm_1(\om)})$=span$\{\phi_1\}$ and $<\phi_1>_{\sscp \lm_1(\om)}\neq 0$. As remarked
in the introduction, since $|\om|=1$, then $\all\leq 1$ and $\all=1$ if and only if $\lm=0$.\\
In particular we can apply Proposition \ref{pr3.1} and continue $\mathcal{G}_{\lm_1(\om)}$ to a
continuous and simple curve without bifurcation points and which locally around any point $s_0>0$
admits a real analytic reparametrization, that is, an injective and continuous map
$\beta_0:(-1,1)\to (s_0-\eps, s_0+\eps)$, $s=\beta(t)$, such that $\beta(0)=s_0$ and
$(\lm(\beta_0(t),\al(\beta_0(t)), \psi(\beta_0(t)))$ is real analytic.
Therefore locally this branch has also the structure of a $1$-dimensional
real analytic manifold and we denote it by,
$$
\mathcal{G}^{(s_1+\dt_1)}=\left\{[0,s_1+\dt_1)\ni s\mapsto (\lm(s),\al(s),\psi(s))\right\},
$$
which satisfies,
$$
\ov{\mathcal{G}^{(s_1+\dt_1)}}=\left\{[0,s_1+\dt_1]\ni s\mapsto (\lm(s),\al(s),\psi(s))\right\},
$$
where, for some $s_1>0$ and $\dt_1>0$, we have:\\
$(A1)_0$ $(\lm(s),\al(s),\psi(s))$ is continuous for $s\in [0,s_1+\dt_1]$;\\
$(A1)_1$ $(\al(s),\psi(s))$ is a positive solution with $\lm=\lm(s)$ for any $s\in [0,s_1+\dt_1]$;\\
$(A1)_2$ $\lm(s)=s$ for $s\leq s_1$, $\lm(s_1)=\lm_1(\om)$ and $\lm(s_*)=\frac1p\Lambda(\om,2p)$
for some $s_*<s_1$;\\
$(A1)_3$ the inclusion $\left\{(\lm,\all,\pl),\,\lm \in [0,\lm_1(\om)] \right\}
\equiv \ov{\mathcal{G}_{\lm_1(\om)}}\subset \mathcal{G}^{(s_1+\dt_1)}$,
holds;\\
$(A1)_4$ $\al(s)\in (0,1),\, \forall\,s\in (0,s_1+\dt_1]$;\\
$(A1)_5$ $\lm(s)>\frac1p\Lambda(\om,2p),\forall\,s \in [s_1,s_1+\dt_1]$;\\
$(A1)_6$ $0\notin \sigma(L_{\ssl(s)}),\, \forall\,s\in (0,s_1+\dt_1)\setminus\{s_1\}$.\\
Clearly $(A1)_5$ holds by $(A1)_2$ and Theorem A. Also, $(A1)_6$ holds since for any $s_0$,
the eigenvalues $\sg_n(s_0)=\sg_n(\al(s_0),\psi(s_0))$ of $L_{\sscp \lm(s_0)}$ are,
after a suitable reparametrization $s=\beta(t)$,
locally real analytic functions (see \cite{but}). Therefore the level sets of a fixed $\sg_n(\beta(t))$
cannot have accumulation points unless $\sg_n(\beta(t))$ is locally constant and consequently,
since $\sg_n(s)$ is (locally up to a reparametrization) real analytic, unless it is constant on $[0,s_1+\dt_1)$. However
we can rule out this case since for $\lm\leq \frac1p\Lambda(\om,2p)$ we have $\sg_1(\all,\pl)>0$ and then we deduce
from $(A1)_2$ that no $\sg_n(s)$ can vanish identically. Therefore it is well defined,
$$
s_2:=\sup\left\{t>s_1\,:\,\al(s)>0\mbox{ and }0\notin \sigma(L_{\ssl(s)}),\;
\forall\,(\lm(s),\al(s),\psi(s))\in \mathcal{G}^{(t)},\,\forall s_1<s<t\right\}.
$$

At this point either $\inf\limits_{s\in(0,s_2)}\al(s)=0$ or $\inf\limits_{s\in(0,s_2)}\al(s)>0$.\\
If $\inf\limits_{s\in(0,s_2)}\al(s)=0$ we are done since we can set $s_\ii=s_2$ and
\beq\label{ii2}
\mathbb{G}_\ii=\mathcal{G}^{(s_\ii)}=\left\{[0,s_\ii)\ni s\mapsto (\lm(s),\al(s),\psi(s))\right\}.
\eeq
If $\inf\limits_{s\in(0,s_2(\om))}\al(s)>0$ the situation is different this time since we don't know
much about the sign of the first eigenvalue for $s>s_1$ and consequently in particular
about the monotonicity of $\all$. However by compactness, Theorem B and since
$\all=1$ implies $\lm=0$, then it is easy to see that necessarily
$\sup\limits_{s\in(0,s_2)}\al(s)<1$. In particular, if $\{\lm_2,\psi_2\}$ denote any limit point of
$\{\lm(s),\psi(s)\}$ as $s\to s_2^{-}$, then by the definition of $s_2$ we must necessarily have that
$Ker(L_{\sscp \lm(s_2)})$ is not empty and then, since $(b)$ of Theorem \ref{thm1.generic} holds,
conclude that in fact $Ker(L_{\sscp \lm_2})=\,$span$\{\phi_2\}$ and
$<\phi_2>_{\sscp \lm(s_2)}\neq 0$. Therefore we can apply once more Proposition \ref{pr3.1} and
conclude that in fact  $(\lm(s),\al(s),\psi(s))$ can be continued as a continuous, locally (up to
reparametrizations) real analytic simple curve. In particular we can argue by induction and for $k\geq 3$ define,
$$
s_k:=\sup\left\{t>s_{k-1}\,:\,\al(s)>0\mbox{ and }0\notin \sg(L_{\ssl(s)}),\;
\forall\,(\lm(s),\al(s),\psi(s))\in \mathcal{G}^{(t)},\,\forall s_{k-1}<s<t\right\}.
$$
If there exists some $k\geq 3$
such that $\inf\limits_{s\in(0,s_k)}\al(s)=0$, then
as in \rife{ii2} we are done.
Otherwise by using Lemma \ref{lemE1}, $(b)$ of Theorem \ref{thm1.generic}, Proposition \ref{pr3.1}
and Theorem A we can find sequences $\{s_k\}$ and $\{\dt_k>0\}$ such that, for any $k\in \mathbb{N}$ we have,
$s_{k+1}>s_{k}>\cdots>s_2>s_1$, $s_{k}+\dt_{k}<s_{k+1}$ and,\\
$(Ak)_0$ $(\lm(s),\al(s),\psi(s))$ is a
continuous and simple curve without bifurcation points (which admits local real analytic
reparamterizations) defined for $s\in [0,s_{k}+\dt_{k}]$ ;\\
$(Ak)_1$ for any $s\in [0,s_k+\dt_k]$, $(\al(s),\psi(s))$ is a positive solution with $\lm=\lm(s)$;\\
$(Ak)_2$ $\lm(s)=s$ for $s<s_1$, $\lm(s_1)=\lm_1(\om)$, $\lm(s_*)=\frac1p\Lambda(\om,2p)$, $s_*< s_1$;\\
$(Ak)_3$ the inclusion
$\left\{(\lm(s),\al(s),\psi(s)),\,s \in [0,s_{k}] \right\}\equiv
\ov{\mathcal{G}^{(s_k)}}\subset \mathcal{G}^{(s_k+\dt_k)}$,
holds;\\
$(Ak)_4$ $\al(s)\in (0,1),\, \forall\,s\in (0,s_k+\dt_k]$;\\
$(Ak)_5$ $\lm(s)>\frac1p\Lambda(\om,2p),\forall\,s \in (s_1,s_k+\dt_k]$;\\
$(Ak)_6$ $0\notin \sg(L_{\ssl(s)}),\, \forall\,s\in (0,s_k+\dt_k)\setminus\{s_1,s_2,\cdots,s_{k}\}$;\\
$(Ak)_7$ $Ker(L_{\sscp \lm(s_k)})=\,$span$\{\phi_k\}$ and $<\phi_k>_{\sscp \lm(s_k)}\neq 0$.\\

Let $s_\ii=\lim\limits_{k\to +\ii}s_k$, we claim that:\\
{\bf Claim:} $\al(s)\to 0^+$ as $s\to s_\ii$.\\
We argue by contradiction and assume that
along an increasing sequence $\{\widehat{s}_j\}$ such that $\widehat{s}_j\to s_{\ii}$,
it holds $\al(\widehat{s}_j)\geq \eps_0>0$ for some $\eps_0>0$. Clearly we can extract a subsequence
$\{s_{k_j}\}\subset \{s_k\}$ such that $s_{k_j}<\widehat{s}_j\leq s_{k_{j+1}}$.
By Theorem \ref{thlambda} we can extract an increasing subsequence of $s$ (which we will not relabel) such that
$\lm(\widehat{s}_{j})\to \widehat{\lm}$
and $\al(\widehat{s}_{j})\to \widehat{\al}\in (0,1]$ as $j\to +\ii$.
Clearly by Lemma \ref{lemE1}, passing to a further subsequence if necessary, we can find
$\widehat{\psi}$ such that $(\widehat{\al},\widehat{\psi})$ is a positive solution for $\lm=\widehat{\lm}$
and $(\al(\widehat{s_{j}}),\psi(\widehat{s_{j}}))\to (\widehat{\al},\widehat{\psi})$ in $C^2_0(\ov{\om})$.
However,
by $(Ak)_5$, $\widehat{\al}$ cannot be $1$ since $\al=1$ if and only if $\lm=0$. Therefore
$\widehat{\al}\in (0,1)$. By Theorem \ref{thm1.generic} we can apply either Lemma \ref{lem1.1} or
Proposition \ref{pr3.1} and conclude that locally around $(\widehat{\al},\widehat{\psi})$
the set of solutions of $\prl$ is a real analytic parametrization of the form
$(\widehat{\lm}(t),\widehat{\al}(t),\widehat{\psi}(t))$, $t\in (-\eps,\eps)$ for some $\eps>0$
with $(\widehat{\lm}(0),\widehat{\al}(0),\widehat{\psi}(0))=(\widehat{\lm},\widehat{\al},\widehat{\psi})$.
In particular for $j$ large enough we can assume without loss of generality that
$(\widehat{\lm}(t),\widehat{\al}(t),\widehat{\psi}(t))$,
$t\in (-\eps,0)$  coincides with $(\lm(s),\al(s),\psi(s))$, $s\in (\widehat{s_{j}},s_{\ii})$.
Let $\{\widehat{\sg}_n\}_{n\in\N}$ be the eigenvalues corresponding to
$(\widehat{\al},\widehat{\psi})$ and $\{\widehat{\sg}_n(t)\}_{n\in\N}$ be those corresponding to
$(\widehat{\al}(t),\widehat{\psi}(t))$.
On one side, since by construction
$0\in \sg(L_{\sscp {\lm(s_{k_j})}})$ and $s_{k_j}<\widehat{s}_j\leq s_{k_{j+1}}$ for any $j$,
then we have that $0\in \sg(L_{\sscp \widehat{\lm}})$.
Indeed, if this was not the case, then, by Lemma \ref{lem1.1} and since the
eigenvalues are isolated, we would have
that there exists a fixed full neighborhood of $0$ with empty intersection with
$\sg(L_{\sscp {\lm(s_{k_j})}})$ for any $j$ large enough, which is a contradiction since the
number of negative eigenvalues is, locally around each positive solution, uniformly bounded.
As a consequence there exists $n\in \N$ such that $\widehat{\sg}_n=0$.
On the other side, since $\widehat{\sg}_n(t)$ is in particular a continuous function of $t$,
by using once more the fact that the eigenvalues are isolated,
possibly passing to a further subsequence if necessary, we must obviously have
$\widehat{\sg}_n(\widehat{t_j})=0$ for some $\widehat{t_j}\to 0^{-}$ as $j\to +\ii$.
Whence $\widehat{\sg}_n$ must vanish identically in $(-\eps,0]$. In particular the $n$-th eigenvalue
of $(\al(s),\psi(s))$ must vanish identically for $s\in (\widehat{s_j},s_{\ii})$ and therefore
in $[0,s_{\ii})$.
This is again a contradiction to $(Ak)_2$
since for $\lm\leq \frac1p\Lambda(\om,2p)$ we have $\sg_1(\all,\pl)>0$ and then
no eigenvalue can vanish identically.
Therefore a contradiction arise which shows that $\al(s)\to 0^+$ as $s\to s_\ii$
and concludes the proof of Lemma \ref{altlem}.\finedim

\bigskip

In view of Lemma \ref{altlem}, to conclude the proof it is enough to show
that, possibly replacing $\mathcal{F}$ of Theorem \ref{thm1.generic} with
a larger but still meager set, then there is a unique value $\lm_\ii(\om,p)$ which can be attained along a subsequence.\\

By a result in \cite{H} (see Example 6.5 and also \cite{ST}), there exists a meager set
$\mathcal{F}^{(1)}_{N,p}\subset \mbox{\rm Diff}^{\,4}(\om_0)$ such that if
$h\in \mbox{\rm Diff}^{\,4}(\om_0)\setminus\mathcal{F}^{(1)}_{N,p}$ and
$\om=h(\om_0)$, then the set of solutions of
$$
-\Delta v=v^p \quad \mbox{in}\;\;\om,\quad  v=0 \quad \mbox{on}\;\;\pa \om
$$
is at most countable. Clearly $\mathcal{F}_{N,p}:=\mathcal{F}\cup \mathcal{F}^{(1)}_{N,p}$ is a meager
set in $\mbox{\rm Diff}^{\,4}(\om_0)$ and
we conclude the proof by showing that if
$h\in \mbox{\rm Diff}^{\,4}(\om_0)\setminus\{\mathcal{F}_N\cup \mathcal{F}^{(1)}_{N,p}\}$ then
there is a unique value $\lm_\ii(\om,p)$ which can be attained along a subsequence.

We argue by
contradiction and assume that
$\frac1p\Lambda(\om,2p)<\lm_{\ii,1}<\lm_{\ii,2}\leq \ov{\lm}(\om,p)$ are attained in the limit along two
distinct subsequences. Since $\al(s)\to 0^+$ as $s\to s_\ii$ and since
$\lm(s)$ is continuous, then it is not difficult to see that
any $\lm\in [\lm_{\ii,1},\lm_{\ii}]$ is attained in the limit
along a subsequence and in particular that there exists a solution of $\prl$ with $\al=0$
for any $\lm\in [\lm_{\ii,1},\lm_{\ii,2}]$. Let $\vl=\bl\pl$, $\bl=\lm^{\frac{p}{p-1}}$ then
\beq\label{LEl}
\graf{-\Delta \vl=\vl^p \quad \mbox{in}\;\;\om\\ \vl=0 \quad \mbox{on}\;\;\pa \om\\
\ino \vl^p=\bl
}
\eeq
We would deduce in particular that, as far as
$h\in \mbox{\rm Diff}^{\,4}(\om_0)\setminus\{\mathcal{F}\cup \mathcal{F}^{(1)}_{N,p}\}$,
the set of solutions of
$$
-\Delta \vl=\vl^p \quad \mbox{in}\;\;\om\quad  \vl=0 \quad \mbox{on}\;\;\pa \om,
$$
is not countable. This is a contradiction since, as mentioned above, if
$h\in \mbox{\rm Diff}^{\,4}(\om_0)\setminus\{\mathcal{F}^{(1)}_{N,p}\}$, then the set of
solution is at most countable.
\finedim

\bigskip
\bigskip

At this point we prove Theorem \ref{lmconvex}. We first state the following straightforward corollary of
some uniqueness results in \cite{DGP}.\\

{\bf Theorem C.} {{\rm (\cite{DGP})}} {\it Let $\om\subset \R^2$ be symmetric and convex with respect to the coordinate
directions $x_i$, $i=1,2$ or $\om=B_{R}\subset \R^{N}$, $N\geq 3$. Then
there exists a unique value of $\lm$, which we denote by
$\lm_*(\om,p)$, such that there exists a solution of {\rm $\prl$} with $\al=0$.
In particular there exists a unique solution $\psi_*$ of {\rm $\prl$} with $\lm=\lm_*(\om,p)$
and $\al=0$.}

\

{\em The proof of Theorem \ref{lmconvex}}\\
By Theorem \ref{spectral1} we see that if either $\om\subset \R^2$ is
convex and symmetric with respect to both directions or  $\om=B_R\subset \R^N$, $N\geq 3$,
and if $0\in \sg(L_{\ssl})$, then $Ker(L_{\ssl})$=span$\{\phi_1\}$ and $<\phi_1>_{\ssl}\neq 0$ always holds.
Therefore for this class of domains the alternative $(i)$-$(ii)$ always holds.
By the uniqueness of the solution for $\al=0$ (see Theorem C) we easily deduce that if $(i)$ holds then
necessarily $\lm_\ii(\om,p)=\lm_*(\om,p)$ and $\psi_{\ii}=\psi_*$.\\
{\em Proof of (j)-(jj).}\\
Both $(j)$ and the left hand side inequality in \rife{030221.1}  are
straightforward consequences of Theorem~1.2 in \cite{BJ3} as stated right after Theorem \ref{lmconvex}. Therefore, we are just left with the proof of
the right hand side inequality in \rife{030221.1}. By using once more Theorem C
 and using the scale invariance of the Rayleigh
quotient in \rife{Sob}, it is not difficult to see that $\psi_*$ is a minimizer of \rife{Sob}
with $t=p+1$. Therefore, it follows from the reverse Holder inequality for minimizers
 (see Theorem 2 in \cite{CRa2}) that
$$
\left(\ino \psi_*^p\right)^2\geq \frac{8\pi}{p+1}
\frac{1}{\Lambda(\om,p+1)}\left(\ino \psi_*^{p+1}\right)^{\frac{2p}{p+1}}
$$
and in view of $2 E_*(\om,p)=\lm_*^p(\om,p)\ino \psi_*^{p+1}$ we readily deduce that
$$
\lm_*(\om,p)\leq \frac{\left(\frac{p+1}{8\pi}\Lambda(\om,p+1)\right)^{\frac{p+1}{2p}}}{2 E_*(\om,p)}
$$
where it is easily seen that if $\om=\mathbb{D}_2$ then the equality holds.
\finedim

\bigskip
\bigskip

\section{\bf On the monotonic continuation of the branch of solutions}
The constrained problem $\prl$ comes with its own variational principle, see section \ref{appD}.
Under a natural connectivity assumption about the set of variational solutions we prove that
$(i)$ of Theorem \ref{alt} holds. Let us define,
$$
\mathbb{S}(\om)=\bigcup_{\lm\in[0,\lm^{**}(\om,p))}\left\{(\lm,\all,\pl)\,:\,(\all,\pl)
\;\mbox{\rm  is a positive variational solution of $\prl$}\right\},
$$
where $\lm^{**}(\om,p)$ has been defined in Theorem A.
We say that ${\bf (H)}$ holds if the following connectivity assumption holds:\\

${\bf (H)}\;$ There exists a subset $\mathbb{S}_{**}(\om)\subseteq \mathbb{S}(\om)$
with the following properties:\\
- for each $\lm\in [0,\lm^{**}(\om,p))$,  $\mathbb{S}_{**}(\om)$ is not empty;\\
- for any pair $\{(\lm_1,\al_1,\psi_{1}) ,(\lm_2,\al_2,\psi_{2})\}\subset \mathbb{S}_{**}(\om)$,
there exists a
continuous map (in the product topology),
$[1,2]\ni t\mapsto  (\lm(t),\al(t),\psi(t))$ such that:\\
$(\lm(1),\al(1),\psi(1))=(\lm_1,\al_1,\psi_{1})$, $(\lm(2),\al(2),\psi(2))=(\lm_2,\al_2,\psi_{2})$ and
$(\lm(t),\al(t),\psi(t))\in \mathbb{S}_{**}(\om)$ for any $t\in [1,2]$.

\bigskip

Then we have,
\bte\label{thm5}
Let $p\in (1,p_N)$ and suppose ${\bf (H)}$ holds for any domain. Then, for most bounded domains $\om\subset \R^N$ of class $C^4$, $|\om|=1$, we have that $\mathbb{S}_{**}(\om)$ is a continuous, simple curve without bifurcation points,
$$
\mathbb{S}_{**}(\om)=\left\{[0,\lm^{**}(\om,p))\ni\lm\mapsto (\all,\pl)\in
(0,1]\times C^{2,r}_{0,+}(\ov{\om}\,)  \right\}
$$
such that
$$
\all\mbox{ is strictly decreasing and } \el\mbox{ is strictly increasing}
$$
and,  as $\lm\nearrow \lm^{**}(\om,p)$, $\all\to 0^+$  and
 $\el\to E_*$. In particular, for $N\geq 3$ it holds $E_*<\dfrac{q}{\lm^{**}(\om,p)}$.
\ete
\proof
We argue as in the very beginning of the proof of Theorem \ref{alt} and deduce
from Lemma~\ref{lemE1} and Lemma \ref{lem1.1} that
we can continue the curve $\mathcal{G}(\om)$ of solutions $(\all,\pl)$ in a right neighborhood
of $\frac1p\Lambda(\om,2p)$ to a larger simple real analytic curve
of solutions with no bifurcation points,
$\mathcal{G}(\om)\subset \mathcal{G}_{\mu}$, $\mu>\frac1p\Lambda(\om,2p)$, such that,  by continuity,
$\all>0$ and $\sg_1(\all,\pl)>0$ for any $\lm<\mu$.
Therefore it is well defined,
$$
\wlm:=\sup\left\{\mu>\frac1p\Lambda(\om,2p)\,:\,\all>0\mbox{ and }\sg_1(\all,\pl)>0,\;
\forall\,(\all,\pl)\in \mathcal{G}_{\mu},\,\forall \lm<\mu\right\}.
$$
Let
\beq\label{Gmu}
\mathbb{G}_{\mu}=\{(\lm,\all,\pl)\,|\,\lm\in[0,\mu), (\all,\pl)\in\mathcal{G}_{\mu}\}.
\eeq
Obviously by Theorem B we have,
\beq\label{GS}
\mathbb{G}_{\frac1p\Lambda(\om,2p)}\subset \mathbb{S}_{**}(\om)
\eeq
and in particular the only solutions for $\lm<\frac1p\Lambda(\om,2p)$ are
the variational solutions. In view of the assumption ${\bf (H)}$, we deduce
the following
\ble\label{lem5.1}
${\bf (H)}$ implies that $\mathbb{G}_{\wlm}\subseteq  \mathbb{S}_{**}(\om)$.
\ele
\proof Suppose by contradiction that there exists $\frac1p\Lambda(\om,2p)\leq \lm_1<\wlm$
such that $(\lm_1,\al_1,\psi_1)\in \mathbb{G}_{\wlm}$ but
$(\lm_1,\al_1,\psi_1)\notin \mathbb{S}_{**}(\om)$. Let us fix $\lm_2\in(\lm_1,\wlm)$ and observe that,
by ${\bf (H)}$, there
exists at least one variational solution
$(\lm_{\sscp \lm_2},\al_{\sscp \lm_2},\psi_{\sscp \lm_2})\in \mathbb{S}_{**}(\om)$. Moreover,
still by ${\bf (H)}$, for fixed
$\lm_0\in [0,\frac1p\Lambda(\om,2p)\,)$, there
exists a continuous path, say $\Gamma$, such that
$$
[0,2]\ni t \xrightarrow {\Gamma} (\lm(t),\al(t),\psi(t))\in
\mathbb{S}_{**}(\om),\mbox{ for any } t\in[0,2],
$$
and
$$
(\lm(2),\al(2),\psi(2))=(\lm_2,\al_{\sscp \lm_2},\psi_{\sscp \lm_2}),\quad
(\lm(0),\al(0), \psi(0))=(\lm_0,\al_{\lm_0},\psi_{\lm_0}),
$$
where $(\al_{\lm_0},\psi_{\lm_0})$ is the unique solution, whence the unique variational solution,
for $\lm=\lm_0$.\\
Because of the gap at $\lm_1$, $(\lm(t),\al(t),\psi(t))\notin \mathbb{G}_{\wlm}$ for those $t$ such that
 $\lm(t)=\lm_1$. We remark that $\lm(t)$ is a continuous function and then it must attain also
 the value $\lm_1$. Therefore the support of
 $\Gamma$ cannot intersect $\mathbb{G}_{\wlm}$ for those $t$ such that $\lm(t)=\lm_1$.
On the other side, since $\mathbb{G}_{\wlm}$ is a real analytic curve with no bifurcation points,
and $\Gamma$ is continuous, then the curve $\Gamma$ cannot intersect $\mathbb{G}_{\wlm}$ at all. Therefore
in particular $(\lm_0,\al_{\lm_0},\psi_{\lm_0})\notin \mathbb{G}_{\wlm}$ which is
a contradiction since the solution at $\lm_0$ is unique and then $(\lm_0,\al_{\lm_0},\psi_{\lm_0})\in \mathbb{G}_{\wlm}$.
\finedim

\bigskip

At this point we deduce the following,

\ble\label{thm1} $\mathbb{G}_{\wlm}$ can be continued to a continuous and piecewise real
analytic simple curve without bifurcation points
$$
\mathbb{G}_{**}=\left\{[0,\lm^{**}(\om,p))\ni\lm\mapsto (\all,\pl)\in
(0,1]\times C^{2,r}_{0,+}(\ov{\om}\,)  \right\}
$$
with the property that $\all$ and $\pl$ are real analytic in $\lm$ with the exception of
at most a strictly increasing set $\{\lm_k\}_{k\in I}$, $I\subseteq \N$, of isolated points, with no finite accumulation points,
excluding possibly  $\lm^{**}(\om,p)$. In particular,
$$
\frac{d\all}{d\lm}<0\mbox{ and }  \frac{dE_{\ssl}}{d\lm}>0, \quad \forall\, 0\leq \lm<\lm^{**}(\om,p),\lm \notin \{\lm_k\}_{k\in \N},
$$
and each $\lm_k$ is an inflection point,
$$
\frac{d\all}{d\lm}\to -\ii\mbox{ and } \frac{d E_{\ssl}}{d\lm}\to +\ii,\;\mbox{ as }\lm \to \lm_k, \,k\in\N.
$$
In particular $\mathbb{S}_{**}(\om)\equiv\mathbb{G}_{**}$.
\ele
\proof   By Lemma \ref{lemE1} and since $\wlm\leq\lm^{**}(\om,p)$,
it is not difficult to see that we can pass to the limit as $\lm\to \wlm$ and conclude that
$\wall=\al_{\sscp \wlm}$, $\wpl=\psi_{\wlm}\in C^{2,r}_{0,+}(\ov{\om}\,)$ and
$\wsg=\sg_1(\wall,\wpl)$ are well defined. In particular,
in view of Lemma \ref{lem5.1}, we see that $(\wall,\wpl)$, being the uniform limit of variational solutions,
is also a variational solution. At this point, if $\wlm=\lm^{**}(\om,p)$, then Lemma \ref{lem1.1} and
Proposition \ref{pr-enrg} yields the desired conclusion.
Therefore we can assume without loss of generality that $\wlm<\lm^{**}(\om,p)$, $\wsg=0$ and $\wall>0$.\\
By Theorem \ref{thm1.generic}
we see that the first eigenvalue $\wsg$ is simple and that, denoting by
$\widehat{\phi}=\phi_{1,\wlm}$ the unique eigenfunction of $\wsg$, it holds
$<\widehat{\phi}>_{\sscp \wlm}\neq 0$. Therefore we can apply Proposition \ref{pr3.1} and conclude that
$\mathcal{G}_{\wlm}$ can be continued to a continuous and piecewise analytic curve with no bifurcation points, say $\mathcal{G}^{(\eps)}$,
defined in a slightly larger interval $[0,\wlm+\eps)$. Please observe that, in the application of Proposition \ref{pr3.1},
we can exclude that $\lm(s)$ is constant, since $\lm$ itself is the free variable of the parametrization for $\lm<\wlm$.
Let $\mathbb{G}_{\sscp \wlm+\eps}$ defined as in \rife{Gmu} just with $\mu=\wlm+\eps$.
In view of ${\bf (H)}$ the same argument of Lemma \ref{lem5.1} shows that
$\mathbb{G}_{\sscp \wlm+\eps}\subseteq \mathbb{S}_{**}(\om)$,
whence, in particular in view of Proposition \ref{prmin}, we have $\sg_1(\all,\pl)\geq 0$ along
$\mathbb{G}_{\sscp \wlm+\eps}$.
At this point we can use \rife{2907.10}
and observe that
$$
\dfrac{\sg_1(s)}{\lm^{'}(s)}=\dfrac{<[\widehat{\phi}]_{\sscp \wlm},\wpl>_{\sscp  \wlm}+\mbox{\rm o}(1)}
{<[\widehat{\phi}]_{\sscp  \wlm}^2>_{\sscp \wlm}+\mbox{\rm o}(1)},\mbox{ as }s\to 0,
$$
where $\widehat{\phi}=\phi_{1,\wlm}$ is the unique eigenfunction of $\wsg$. We can assume without loss of generality that $\lm^{'}(s)>0$ for $s<0$.
By using once more $<\widehat{\phi}>_{\sscp \wlm}\neq 0$ and \rife{2907.5}, then we conclude that $<[\widehat{\phi}]_{\sscp \wlm},\wpl>_{\sscp  \wlm}\neq 0$
and actually, since $\sg_1(s)> 0$ and $\lm^{'}(s)>0$ for $s<0$, that necessarily $<[\widehat{\phi}]_{\sscp \wlm},\wpl>_{\sscp  \wlm}> 0$.
As a consequence, since $\sg_1(s)\geq 0$ for $s>0$ small enough, then in particular $\sg_1(s)>0$ and consequently $\lm^{'}(s)>0$ for any $s\neq 0$ small enough. Therefore
$\mathcal{G}^{(\eps)}\equiv \mathcal{G}_{\sscp \wlm+\eps}$ is a continuous, piecewise
real analytic curve of solutions of $\prl$ with $\lm\in [0,\wlm+\eps)$ and $\sg_1(\pl,\all)>0$ in $[0,\wlm+\eps)\setminus\{\wlm\}$. In particular, in view of
\rife{2907.1} and since $\lm^{'}(s)\to 0^+$ as $s\to 0$ we find that,
$$
\frac{d\all}{d\lm}=\frac{\al^{'}(s)}{\lm^{'}(s)}=\frac{-\wlm<\wpl>_{\sscp \wlm}+\mbox{o}(1)}{\lm^{'}(s)}\to -\ii,\mbox{ as }s\to 0.
$$
On the other side we have,
$$
\dfrac{\ino \rh_{\sscp \wlm} \widehat{\phi}}{\ino (\rh_{\sscp \wlm})^{\frac1q}}=
\wall<\widehat{\phi}>_{\sscp \wlm}+\wlm<\wpl \,[\widehat{\phi}]_{\sscp \wlm}\!\!>_{\sscp \wlm}+
\wlm<\wpl>_{\sscp \wlm}<\widehat{\phi}>_{\sscp \wlm},
$$
which, by using once more \rife{2907.5}, shows that $\ino \rh_{\sscp \wlm}\widehat{\phi}>0$. Therefore we conclude that,
$$
\frac{d E_{\ssl}}{d\lm}=\ino \rh_{\lm} \frac{d\pl}{d\lm}=\ino \rh_{\sscp \lm(s)}\frac{\psi^{'}(s)}{\lm^{'}(s)}=
\ino \left(\rh_{\sscp \wlm}+\mbox{o}(1)\right)\frac{\widehat{\phi}+\mbox{o}(1)}{\lm^{'}(s)}\to +\ii,\mbox{ as }s\to 0.
$$

At this point we can simply iterate the argument above. Indeed, by defining

$$
\widehat{\lm}_2=\widehat{\lm}_2(\om):=\sup\{\widehat{\lm}_1<\mu<\lm^{**}(\om,p)\,:\,\all>0\mbox{ and }\sg_1(\pl,\all)>0,\;
\forall\,(\all,\pl)\in \mathcal{G}_{\mu}\},
$$
we see that we are left with the same two possibilities, that is, either $\widehat{\lm}_2=\lm^{**}(\om,p)$, and we would be done, or
$\widehat{\lm}_2<\lm^{**}(\om,p)$ and we could repeat the continuation argument. At this point,
the same inductive argument adopted in the proof of Theorem \ref{alt} shows that
$\mathbb{G}_{\sscp \wlm+\eps}$ can be continued in the full interval
$(0,\lm^{**}(\om,p))$ to a
continuous curve $\mathbb{G}_{\sscp \wlm+\eps}$ $(\lm,\all,\pl)$, where $(\all,\pl)$ are
positive variational solutions of $\prl$ with the property that $\all$ and $\pl$ are
real analytic in $\lm$ with the exception of
at most a strictly increasing sequence $\{\lm_k\}_{k\in\N}$ of isolated points such that
$\lm_k\to \lm^{**}(\om,p)$.
Indeed if an accumulation point of $\lm_k$ would exists which is not $\lm^{**}(\om,p)$,
then we would obtain a contradiction to the real analiticity of the parametrization
around positive solutions. This is done exactly as in the proof Theorem \ref{alt} and we skip
this part here to avoid repetitions.
In particular we have,
$$
\frac{d\all}{d\lm}<0\mbox{ and }  \frac{dE_{\ssl}}{d\lm}>0, \quad \forall \lm \neq \lm_k,\,k\in\N,
$$

and each $\lm_k$ is an inflection point,

$$
\frac{d\all}{d\lm}\to -\ii\mbox{ and } \frac{d E_{\ssl}}{d\lm}\to +\ii,\;\mbox{ as }\lm \to \lm_k, \,k\in\N.
$$
\finedim

\bigskip
\bigskip
Putting $\mathbb{G}_{**}=\left.\mathbb{G}_{\mu}\right|_{\lm^{**}(\om,p)}$, we see that
Lemma \ref{thm1} essentially concludes the proof. Indeed, the bound for $\el$ is an
immediate consequence of \rife{enrgbound} below and we just miss the fact that
$\mathbb{G}_{**}\equiv \mathbb{S}_{**}(\om)$. However this is easy to check, since by construction
$\mathbb{G}_{**}\subseteq \mathbb{S}_{**}(\om)$, while if
$\mathbb{S}_{**}(\om)\setminus \mathbb{G}_{**}$ were not empty, then by ${\bf (H)}$
we could always find a continuous
path connecting points on {$\mathbb{S}_{**}(\om)\setminus \mathbb{G}_{**}$} whose $\lm$ satisfies
$\frac1p\Lambda(\om,2p)\leq \lm<\lm^{**}(\om,p)$ with
any one of those in $\mathbb{S}_{**}(\om)\cap \mathbb{G}_{**}$ with $\lm_0<\frac1p\Lambda(\om,2p)$.
As in Lemma \ref{lem5.1} this would cause $\mathbb{G}_{**}$
to have a bifurcation point at some $(\lm_0,\al_{\lm_0},\psi_{\lm_0})<\frac1p\Lambda(\om,2p)$,
which is impossible.
\finedim

\bigskip
\bigskip

\section{\bf Applications: a nonlinear eigenvalue problem}
We present here an application to a classic nonlinear eigenvalue problem, which was actually our starting motivation to pursue this approach. In a sense, we obtain a ``desingularization" of the blow up limit of problem \rife{ul} below.
For $p\in(1,p_N)$, let us consider classical solutions of
\beq\label{ul}
\graf{-\Delta u=\mu (1+u)^p \quad \mbox{in}\;\;\om\\ u=0 \quad \mbox{on}\;\;\pa \om}.
\eeq
By assuming $p<\frac{N+2}{N-2}$ and $\um$ to be
a mountain pass solution of \rife{ul}, it was proved in \cite{GM} that along a subsequence and
as $\mu\to 0^+$ then
$\mu^{\frac{1}{p-1}}\um \to U_{\ii}$ where $U_\ii>0$ in $\om$ is a mountain pass solution of
\beq\label{Uinfty}
-\Delta U_\ii=U_\ii^p, \mbox{ in }\om, \quad U_{\ii}=0 \mbox{ on }\pa \om.
\eeq

Here we deduce the following generic
result about the Rabinowitz (\cite{Rab}) unbounded continuum of solutions of \rife{ul}.

\bte\label{Rab}
Let $p\in(1,p_N)$, $\om\subset \R^N$ be as in Theorem \ref{alt} and let $(\lm(s),\all(s),\pl(s))$, $s\in [0,s_\ii)$,
denote the
curve of positive solutions of {\rm $\prl$} satisfying either $(i)$ or $(ii)$ in Theorem~\ref{alt}.
Then, there exists $\lm_\ii(\om,p)\in (\frac1p\Lambda(\om,2p),\ov{\lm}(\om,p)]$ such that the Rabinowitz
unbounded continuum of solutions of {\rm \rife{ul}} is a continuous parametrization
$$
\mathcal{R}_{\ii}=\{[0,s_\ii)\ni s \mapsto (\mu(s),u(s))\in [0,+\ii)\times C^{2,r}_0(\ov{\om}\,)\},
$$
where
$$(\mu(s),u(s))=(\lm(s)\al^{p-1}(s),\frac{\lm(s)}{\al(s)}\psi(s)),
$$
and $\al(s)\in (0,1]$, $\lm(s)\in (\frac1p\Lambda(\om,2p),\ov{\lm}(\om,p))$,
$(\mu(0),u(0))=(0,0)$.\\
Moreover, as $s\to s_\ii$, $\al(s)\to 0^+$,
$\lm(s)\to\lm_\ii(\om,p)$, so that
$\mu(s)\to 0^+$ and in particular for any sequence $s_n\to s_\ii$
there exists a subsequence $\{t_n\}\subseteq \{s_n\}$,
such that $\lm^{\frac{1}{p-1}}(t_n)\al(t_n)u(t_n)\to U_{\ii}$ in $C^{2}_0(\ov{\om}\,)$,
where $U_\ii$ is a solution of {\rm \rife{Uinfty}}.
\ete

\smallskip

\brm\label{remb2}{\it  Obviously, in view of Theorem \ref{lmconvex},
if $\om\subset \R^2$ is convex and symmetric with respect to the coordinate directions or either
$\om=\mathbb{D}_N$, $N\geq 3$,
then, as $s\to s_\ii$, $\lm^{\frac{1}{p-1}}(s)\al(s)u(s)\to U_{*}$ where $U_{*}$ is the unique
solution of {\rm \rife{Uinfty}}. However, in these particular cases much more is known about
$\mathcal{R}_{\ii}$, see for example \cite{HK}, \cite{Kor}. Also, it is well known that $\mu\to 0^+$
along the branch of non minimal solutions {\rm(\cite{Lions})}.\\
The interest of Theorem \cite{Rab} relies in the "desingularization" of
the limit $\mu\to 0^+$ for non minimal solutions, that is, the fact that,
as a corollary of Theorem \ref{alt}, in a generic sense along $\mathcal{R}_{\ii}$ any non minimal unbounded
sequence of solutions must be asymptotically "almost" proportional to a solution of {\rm \rife{Uinfty}}.
}
\erm

{\it Proof of Theorem {\rm \ref{Rab}}.} We first point out that in case $(i)$ of Theorem \ref{alt} we just have $\lm(s)=s$, $\al(s)=\al_{s}$, $\psi(s)=\psi_s$. Now, it is readily seen that for any $s$, $u({s})=\frac{\lm(s)}{\al(s)}\psi(s)$ is a solution of
\rife{ul} with $\mu=\mu(s)=\lm(s)\al^{p-1}(s)$. Since by Theorem \ref{alt} $(\lm(s),\all(s),\pl(s))$ is continuous,
we deduce that
$$
[0,s_\ii)\ni s\mapsto (\mu(s),u(s))=\left(\lm(s)\al^{p-1}(s),\frac{\lm(s)}{\al(s)}\psi(s)\right)
$$
is a continuous (not necessarily simple) curve. Since in particular $\al(0)=1$, $\lm(0)=0$, $\psi(0)=G[1]$,
we are just left to prove the assertion about the limit $s\to s_\ii$. By Theorem \ref{alt}
we have $\al(s)\to 0^+$ and $\lm(s)\to \lm_{\ii}(\om,p)$ as $s\to s_\ii$
and for any sequence $s_n\to s_\ii$ there exists a subsequence $\{t_n\}\subseteq \{s_n\}$
such that
$\psi(t_n)\to \psi_{\ii}$ in $C_0^2(\ov{\om}\,)$, for some $\psi_{\ii}$ which solves
{\rm $\prl$} with $\lm=\lm_\ii(\om,p)$ and $\al=0$. Clearly
$U(s)=\lm^{\frac{1}{p-1}}(s)\al(s)u_s=\lm^{\frac{p}{p-1}}(s)\psi(s)$ is a solution of
$$
\graf{-\Delta U(s)=\lm^{\frac{p}{p-1}}(s)(\al(s)+\lm^{-\frac{1}{p-1}}(s) U(s))^p \quad \mbox{in}\;\;\om
\\ U(s)=0 \quad \mbox{on}\;\;\pa \om}
$$
and the convergence of $U(t_n)$ follows by standard elliptic estimates.
\finedim

\bigskip
\bigskip

\appendix

\section{Variational solutions}\label{appD}
Problem  $\prl$ arises as the Euler-Lagrange equation of the constrained
minimization principle {\bf (VP)} below for the plasma densities {\rm $\rh\in L^{1}(\om)$},
which, for $p>1$, is equivalent to the variational formulation of \fbi.
We shortly discuss here the variational solutions of \fbi\, and $\prl$ and their equivalence and
refer to \cite{BeBr} for a detailed discussion of this point. In this context $\all$ is the Lagrange multiplier related to the "mass" constraint $\ino \rl=1$ while the Dirichlet energy is the density interaction energy,
$$
\mathcal{E}(\rh)=\frac12\ino \rh G[\rh],
$$
which is easily seen to coincide with $\el$ whenever $\pl=G[\rl]$, that is, $\el=\mathcal{E}(\rl)$.\\
For any
$$
\rh\in\mathcal{P}_{\sscp \om}:=\left\{\rh\in L^{1+\frac1p}(\om)\,|\,\rh\geq 0\;\mbox{a.e. in}\;\om \right\},
$$
and $\lm\geq  0$, we define the free energy,
\beq\label{jeil}
J_{\ssl}(\rh)=
{\scriptstyle \frac{p}{p+1}}\ino (\rh)^{1+\frac{1}{p}}-\frac\lm 2 \ino \rho G[\rho].
\eeq
Let us consider the variational principle,
$$
\mathcal{J}(\lm)=\inf\left\{ J_{\ssl}(\rh)\,:\,\rh\in \mathcal{P}_{\sscp \om}, \ino \rh=1\right\}.\qquad \qquad \mbox{\bf (VP)}
$$
It has been shown in \cite{BeBr, Te2} that for each $\lm>0$ there exists at least one $\rl$ which solves the {\bf (VP)}.
In particular, (\cite{BeBr}) $\al=\all\in\R$ arises as the Lagrange multiplier relative to the constraint $\ino \rl=1$ and if
$\all\geq 0$, then
any minimizer $\rl$ yields a solution $(\all,\pl)$ of $\prl$ where $\pl=G[\rl]$.
Any such solution is called a variational solution of $\prl$.\\
Solutions of \fbi\, are also found in \cite{BeBr} as minimizers of the
functional $\Psi_{I}(v)$, $v\in\mathcal{H}_I$ defined as follows
\beq\label{psiI}
\Psi_{I}(v)=\frac12\ino |\nabla v|^2-\frac{1}{p+1}\left(\ino (v)_+^{p+1}\right)+ Iv(\pa \om),
\eeq
on
$$
\mathcal{H}_I=\left\{v\in H\,|\,\ino (v)_+^p=I\right\},
$$
where $I>0$ and $H$ is the space of $H^{1}(\om)$ functions whose boundary trace is constant.
For fixed $I>0$, a variational solution of \fbi\, is a solution of \fbi\, which is also a minimizer of
$\Psi_I$ on $\mathcal{H}_I$. It has been shown in \cite{BeBr} that, for $p\in(1,p_N)$, at least
one variational solution exists for each $I>0$. The proof of Theorem A in the introduction can be found in
Appendix A of \cite{BJ1}.\\

Finally we show that for positive variational solutions of {\rm $\prl$} and $N\geq 3$ the following holds:
\beq\label{enrgbound}
E_{\ssl}\leq \frac{q}{\lm}(1-\all),\,\forall \lm>0.
\eeq
Indeed, putting $c_p=\frac{p}{p+1}$, we see that,
$$
\mathcal{J}(\lm)=-1+c_p\ino (\rl)^{1+\frac1p}-\lm E_{\ssl}=-1+c_p\ino \rl (\all+\lm\pl)-\lm E_{\ssl}=
$$
\beq\label{1806.1}
-1+c_p\all+c_p2\lm E_{\ssl}-\lm E_{\ssl}=-1+c_p\all +\frac{p-1}{p+1}\lm E_{\ssl}.
\eeq
Since $\mathcal{J}(\lm)$ is decreasing, then we also find that,
$$
-1+c_p=\mathcal{J}(0)\geq \mathcal{J}(\lm)=-1+c_p\all +\frac{p-1}{p+1}\lm E_{\ssl},
$$
which immediately yields \rife{enrgbound}.

\end{document}